\input amstex
\documentstyle {amsppt}
\pagewidth{12.5 cm} \pageheight{19 cm} \magnification \magstep1
\topmatter
\title Total positivity in the \\ De Concini-Procesi Compactification \endtitle
\author Xuhua He\endauthor

\address Department of Mathematics, M.I.T., Cambridge, MA 02139\endaddress

\email hugo\@math.mit.edu \endemail

\thanks 2000 {\it Mathematics Subject Classification.} Primary 20G20; Secondary
14M15. \endthanks

\abstract We study the nonnegative part $\overline{G_{>0}}$ of the
De Concini-Procesi compactification of a semisimple algebraic
group $G$, as defined by Lusztig. Using positivity properties of
the canonical basis and parametrization of flag varieties, we will
give an explicit description of $\overline{G_{>0}}$. This answers
the question of Lusztig in [L4]. We will also prove that
$\overline{G_{>0}}$ has a cell decomposition which was conjectured
by Lusztig.
\endabstract
\endtopmatter
\document

\define\opp{\bowtie}
\define\po{\text{\rm pos}}

\define\supp{\text{\rm supp}}
\define\End{\text{\rm End}}
\define\diag{\text{\rm diag}}
\define\Lie{\text{\rm Lie}}

\define\Ad{\text{\rm Ad}}
\redefine\i{^{-1}}
\redefine\ge{\geqslant}
\redefine\le{\leqslant}

\define\cb{\Cal B}

\define\cp{\Cal P}

\define\car{\Cal R}

\define\cv{\Cal V}

\define\a{\alpha}

\head 0. Introduction \endhead

Let $G$ be a connected split semisimple algebraic group of adjoint
type over $\bold R$. We identify $G$ with the group of its $\bold
R$-points. In [DP], De Concini and Procesi defined a
compactification $\bar{G}$ of $G$ and decomposed it into strata
indexed by the subsets of a finite set $I$. We will denote these
strata by $\{Z_J \mid J \subset I\}$. Let $G_{>0}$ be the set of
strictly totally positive elements of $G$ and $G_{\ge 0}$ be the
set of totally positive elements of $G$ (see [L1]). We denote by
$\overline{G_{>0}}$ the closure of $G_{>0}$ in $\bar{G}$. The main
goal of this paper is to give an explicit description of
$\overline{G_{>0}}$ (see 3.14). This answers the question in [L4,
9.4]. As a consequence, I will prove in 3.17 that
$\overline{G_{>0}}$ has a cell decomposition which was conjectured
by Lusztig.

To achieve our goal, it is enough to understand the intersection
of $\overline{G_{>0}}$ with each stratum. We set $Z_{J, \ge
0}=\overline{G_{>0}} \bigcap Z_J$. Note that $Z_I=G$ and $Z_{I,
\ge 0}=G_{\ge 0}$. We define $Z_{J, >0}$ as a certain subset of
$Z_{J, \ge 0}$ analogous to $G_{>0}$ for $G_{\ge 0}$ (see 2.6).
When $G$ is simply-laced, we will prove in 2.7 a criterion for
$Z_{J, >0}$ in terms of its image in certain representations of
$G$, which is analogous to the criterion for $G_{>0}$ in [L4,
5.4]. As Lusztig pointed out in [L2], although the definition of
total positivity was elementary, many of the properties were
proved in a non-elementary way, using canonical bases and their
positivity properties. Our theorem 2.7 is an example of this
phenomenon. As a consequence, we will see in 2.9 that $Z_{J, \ge
0}$ is the closure of $Z_{J, >0}$ in $Z_J$.

Note that $Z_J$ is a fiber bundle over the product of two flag
manifolds. Then understanding $Z_{J, \ge 0}$ is equivalent to
understanding the intersection of $Z_{J, \ge 0}$ with each fiber.
In 3.5, we will give a characterization of $Z_{J, \ge 0}$ which is
analogous to the elementary fact that $G_{\ge 0}=\bigcap\limits_{g
\in G_{>0}} g\i G_{>0}$. It allows us to reduce our problem to the
problem of understanding certain subsets of some unipotent groups.
Using the parametrization of the totally positive part of the flag
varieties (see [MR]), we will give an explicit description of the
subsets of $G$(see 3.7). Thus our main theorem can be proved.

\head 1. Preliminaries \endhead

\subhead 1.1 \endsubhead We will often identify a real algebraic
variety with the set of its $\bold R$-rational points. Let $G$ be
a connected semisimple adjoint algebraic group defined and split
over $\bold R$, with a fixed \'{e}pinglage $(T, B^+, B^-, x_i,
y_i; i \in I)$ (see [L1, 1.1]). Let $U^+, U^-$ be the unipotent
radicals of $B^+, B^-$. Let $X$ (resp. $Y$) be the free abelian
group of all homomorphism of algebraic groups $T@>>>\bold R^*$
(resp. $\bold R^*@>>>T$) and $< , >: Y \times X@>>>\bold Z$ be the
standard pairing. We write the operation in these groups as
addition. For $i \in I$, let $\alpha_i \in X$ be the simple root
such that $t x_i(a) t \i=x_i(a)^{\alpha_i(t)}$ for all $a \in
\bold R, t \in T$ and $\alpha_i ^\vee \in Y$ be the simple coroot
corresponding to $\alpha_i$. For any root $\alpha$, we denote by
$U_{\alpha}$ the root subgroup corresponding to $\alpha$.

There is a unique isomorphisms $\psi : G@>\sim>>G^{\text{opp}}$
(the opposite group structure) such that

$\psi\bigl(x_i(a)\bigr)=y_i(a)$, $\psi\bigl(y_i(a)\bigr)=x_i(a)$
for all $i\in I$, $a \in \bold R$ and $\psi(t)=t$, for all $t \in
T$.

If $P$ is a subgroup of $G$ and $g \in G$, we write $^gP$ instead
of $gPg\i$.

For any algebraic group $H$, we denote the Lie algebra of $H$ by
$\Lie(H)$ and the center of $H$ by $Z(H)$.

For any variety $X$ and an automorphism $\sigma$ of $X$, we denote
the fixed point set of $\sigma$ on $X$ by $X^{\sigma}$.

For any group, We will write 1 for the identity element of the
group.

For any finite set $X$, we will write $|X|$ for the cardinal of
$X$.

\subhead 1.2 \endsubhead Let $N(T)$ be the normalizer of $T$ in
$G$ and ${\dot{s_i}=x_i(-1)y_i(1)x_i(-1) \in N(T)}$ for $i \in I$.
Set $W=N(T)/T$ and $s_i$ to be the image of $\dot{s_i}$ in $W$.
Then $W$ together with $(s_i)_{i \in I}$ is a Coxeter group.

Define an expression for $w \in W$ to be a sequence $\hbox{\bf
w}=(w_{(0)}, w_{(1)}, \ldots, w_{(n)})$ in $W$, such that
$w_{(0)}=1$, $w_{(n)}=w$ and for any $j=1, 2, \ldots, n$,
$w_{(j-1)} \i w_{(j)}=1 $ or $s_i$ for some $i \in I$. An
expression $\hbox{\bf w}=(w_{(0)}, w_{(1)}, \ldots, w_{(n)})$ is
called reduced if $w_{(j-1)}<w_{(j)}$ for all $j=1, 2, \ldots, n$.
In this case, we will set $l(w)=n$. It is known that $l(w)$ is
independent of the choice of the reduced expression. Note that if
$\hbox{\bf w}$ is a reduced expression of $w$, then for all $j=1,
2, \ldots, n$, $w_{(j-1)}\i w_{(j)}=s_{i_j}$ for some $i_j \in I$.
Sometimes we will simply say that $s_{i_1} s_{i_2} \cdots s_{i_n}$
is a reduced expression of $w$.

For $w \in W$, set
$\dot{w}=\dot{s_{i_1}}\dot{s_{i_1}}\cdots\dot{s_{i_n}}$ where
$s_{i_1} s_{i_2} \cdots s_{i_n}$ is a reduced expression of $w$.
It is well known that $\dot{w}$ is independent of the choice of
the reduced expression $s_{i_1}s_{i_2} \cdots s_{i_n}$ of $w$.

Assume that $\hbox{\bf w}=(w_{(0)}, w_{(1)}, \ldots, w_{(n)})$ is
a reduced expression of $w$ and $w_{(j)}=w_{(j-1)} s_{i_j}$ for
all $j=1, 2, \ldots, n$. Suppose that $v \le w$ for the standard
partial order in $W$. Then there is a unique sequence $\hbox{\bf
v}_+=(v_{(0)},v_{(1)}, \ldots, v_{(n)})$ such that $v_{(0)}=1,
v_{(n)}=v, v_{(j)} \in \{v_{(j-1)},\, v_{(j-1)}s_{i_j}\}$ and
$v_{(j-1)}< v_{(j-1)}s_{i_j}$ for all $j=1, 2, \ldots, n$ (see
[MR, 3.5]). $\hbox{\bf v}_+$ is called the positive subexpression
of $\hbox{\bf w}$. We define $$\eqalignno{J^+_{\hbox{\bf v}_+}
&=\{j \in \{1, 2, \ldots, n\} \mid v_{(j-1)}<v_{(j)}\}, \cr
J^{\circ}_{\hbox{\bf v}_+} &=\{j \in \{1, 2, \ldots, n\} \mid
v_{(j-1)}=v_{(j)}\}.}$$

Then by the definition of $\hbox{\bf v}_+$, we have $\{1, 2,
\ldots, n\}=J^+_{\hbox{\bf v}_+} \sqcup J^{\circ}_{\hbox{\bf
v}_+}$.

\subhead 1.3 \endsubhead Let $\cb$ be the variety of all Borel
subgroups of $G$. For $B, B'$ in $\cb$, there is a unique $w \in
W$, such that $(B, B')$ is in the $G$-orbit on $\cb \times \cb$
(diagonal action) that contains $(B^+, ^{\dot{w}}B^+)$. Then we
write $\po(B, B')=w$. By the definition of $\po$, $\po(B,
B')=\po(^g B, ^g B')$ for any $B, B' \in \cb$ and $g \in G$.

For any subset $J$ of $I$, let $W_J$ be the subgroup of $W$
generated by $\{s_j \mid j \in J\}$ and $w^J_0$ be the unique
element of maximal length in $W_J$. (We will simply write $w^I_0$
as $w_0$.) We denote by $P_J$ the subgroup of $G$ generated by
$B^+$ and by $\{y_j(a) \mid j \in J, a \in \bold R\}$ and denote
by $\cp^J$ the variety of all parabolic subgroups of $G$
conjugated to $P_J$. It is easy to see that for any parabolic
subgroup $P$, $P \in \cp^J$ if and only if $\{\po(B_1, B_2) \mid
B_1, B_2 \hbox{ are Borel subgroups of } P \}=W_J$.

\subhead 1.4 \endsubhead For any parabolic subgroup $P$ of $G$,
define $U_P$ to be the unipotent radical of $P$ and $H_P$ to be
the inverse image of the connected center of $P/U_P$ under
$P@>>>P/U_P$. If $B$ is a Borel subgroup of $G$, then so is
$$P^B=(P \cap B) U_P.$$

It is easy to see that for any $g \in H_P$, we have $^g
(P^B)=P^B$. Moreover, $P^B$ is the unique Borel subgroup $B'$ in
$P$ such that $\po(B, B') \in W^J$, where $W^J$ is the set of
minimal length coset representatives of $W/W_J$ (see [L5,
3.2(a)]).

Let $P, Q$ be parabolic subgroups of $G$. We say that $P, Q$ are
opposed if their intersection is a common Levi of $P, Q$. (We then
write $P \opp Q$.) It is easy to see that if $P \opp Q$, then for
any Borel subgroup $B$ of $P$ and $B'$ of $Q$, we have $\po(B, B')
\in W_J w_0$.

For any subset $J$ of $I$, define $J^* \subset I$ by $\{Q \mid Q
\opp P \hbox{ for some } P \in \cp^J\}=\cp^{J^*}$. Then we have
$(J^*)^*=J$. Let $Q_J$ be the subgroup of $G$ generated by $B^-$
and by $\{x_j(a) \mid j \in J, a \in \bold R\}$. We have $Q_J \in
\cp^{J^*}$ and $P_J \opp Q_J$. Moreover, for any $P \in \cp^J$, we
have $P=^g P_J$ for some $g \in G$. Thus $\psi(P)=^{\psi(g) \i}
Q_J \in \cp^{J^*}$.

\subhead 1.5 \endsubhead Recall the following definitions from
[L1].

For any $w \in W$, assume that $w=s_{i_1} s_{i_2} \cdots s_{i_n}$
is a reduced expression of $w$. Define $\phi^{\pm}: R^n_{\ge
0}@>>>U^{\pm}$ by
$$\eqalign{\phi^+(a_1, a_2, \ldots, a_n)=x_{i_1}(a_1) x_{i_2}(a_2) \cdots
x_{i_n}(a_n), \cr \phi^-(a_1, a_2, \ldots, a_n)=y_{i_1}(a_1)
y_{i_2}(a_2) \cdots y_{i_n}(a_n). \cr}$$

Let $U^{\pm}_{w, \ge 0}=\phi^{\pm}(R^n_{\ge 0}) \subset U^{\pm}$,
$U^{\pm}_{w, >0}=\phi^{\pm}(R^n_{>0}) \subset U^{\pm}$. Then
$U^{\pm}_{w, \ge 0}$ and $U^{\pm}_{w, >0}$ are independent of the
choice of the reduced expression of $w$. We will simply write
$U^{\pm}_{w_0, \ge 0}$ as $U^{\pm}_{\ge 0}$ and $U^{\pm}_{w_0,
>0}$ as $U^{\pm}_{>0}$.

$T_{>0}$ is the submonoid of $T$ generated by the elements
$\chi(a)$ for $\chi \in Y$ and $a \in \bold R_{>0}$.

$G_{\ge 0}$ is the submonoid $U^+_{\ge 0} T_{>0} U^-_{\ge
0}=U^-_{\ge 0} T_{>0} U^+_{\ge 0}$ of $G$.

$G_{>0}$ is the submonoid $U^+_{>0} T_{>0} U^-_{>0}=U^-_{>0}
T_{>0} U^+_{>0}$ of $G_{\ge 0}$.

$\cb_{>0}$ is the subset $\{^u B^- \mid u \in U^+_{>0} \}=\{^u B^+
\mid u \in U^-_{>0} \}$ of $\cb$ and $\cb_{\ge 0}$ is the closure
of $\cb_{>0}$ in the manifold $\cb$.

For any subset $J$ of $I$, $\cp^J_{>0}=\{P \in \cp^J \mid \exists
B \in \cb_{>0}, \text{ such that } B \subset P\}$ and $\cp^J_{\ge
0}=\{P \in \cp^J \mid \exists B \in \cb_{\ge 0}, \text{ such that
} B \subset P\}$ are subsets of $\cp^J$.

\subhead 1.6 \endsubhead For any $w, w' \in W$, define
$$\car_{w, w'}=\{B \in \cb \mid \po(B^+, B)=w', \po(B^-, B)=w_0
w \}.$$ It is known that $\car_{w, w'}$ is nonempty if and only if
$w \leqslant w'$ for the standard partial order in $W$(see [KL]).
Now set $$\car_{w, w', >0}=\cb_{\geqslant 0} \cap \car_{w, w'}.$$
Then $\car_{w, w', >0}$ is a connected component of $\car_{w, w'}$
and is a semi-algebraic cell(see [R2, 2.8]). Furthermore,
$\cb=\bigsqcup\limits_{w \leqslant w'} \car_{w, w'}$ and
$\cb_{\geqslant 0}=\bigsqcup\limits_{w \leqslant w'} \car_{w, w',
>0}$. Moreover, for any $u \in U^+_{w \i, >0}$, we have $^u
\car_{w, w', >0} \subset \car_{1, w', >0}$ (see [R2, 2.2]).

Let $J$ be a subset of $I$. Define $\pi^J: \cb@>>>\cp^J$ to be the
map which sends a Borel subgroup to the unique parabolic subgroup
in $\cp^J$ that contains the Borel subgroup. For any $w, w' \in W$
such that $w \leqslant w'$ and $w' \in W^J$, set $\cp^J_{w,
w'}=\pi^J(\car_{w, w'})$ and $\cp^J_{w, w', >0}=\pi^J(\car_{w, w',
>0})$. We have $\cp^J_{\geqslant 0}=\bigsqcup\limits_{w \leqslant
w', w' \in W^J} \cp^J_{w, w', >0}$ and $\pi^J \mid _{\car_{w, w',
>0}}$ maps $\car_{w, w', >0}$ bijectively onto $\cp^J_{w, w', >0}$
(see [R1, chapter 4, 3.2]). Hence, for any $u \in U^+_{w \i, >0}$,
we have $^u \cp^J_{w, w', >0}= \pi^J(^u \car_{w, w', >0}) \subset
\pi^J(\cp^J_{1, w', >0})$.

\subhead 1.7 \endsubhead Define $\pi_T: B^- B^+@>>>T$ by $\pi_T(u
t u')=t$ for $u \in U^-, t \in T, u' \in U^+$. Then for $b_1 \in
B^-, b_2 \in B^- B^+, b_3 \in B^+$, we have $\pi_T(b_1 b_2
b_3)=\pi_T(b_1) \pi_T(b_2) \pi_T(b_3)$.

Let $J$ be a subset of $I$. We denote by $\Phi^+_J$ the set of
roots that are linear combination of $\{\alpha_j \mid j \in J\}$
with nonnegative coefficients. We will simply write $\Phi^+_I$ as
$\Phi^+$ and we will call a root $\alpha$ positive if $\alpha \in
\Phi^+$. In this case, we will simply write $\alpha>0$. Define
$U^+_J$ to be the subgroup of $U^+$ generated by $\{U_{\alpha}
\mid \alpha \in \Phi^+_J\}$ and $'U^+_J$ to be the subgroup of
$U^+$ generated by $\{U_{\alpha} \mid \alpha \in
\Phi^+-\Phi^+_J\}$. Then $U^- \times T \times 'U^+_J \times U^+_J$
is isomorphic to $B^- B^+$ via $(u, t, u_1, u_2) \mapsto u t u_1
u_2$. Now define $\pi_{U^+_J}: B^- B^+@>>>U^+_J$ by $\pi_{U^+_J}(u
t u_1 u_2)=u_2$ for $u \in U^-, t \in T, u_1 \in 'U^+_J$ and $u_2
\in U^+_J$. (We will simply write $\pi_{U^+_I}$ as $\pi_{U^+}$.)
Note that $U^- T \cdot U^- T 'U^+_J=U^- T 'U^+_J$. Thus it is easy
to see that for any $a, b \in G$ such that $a, ab \in B^- B^+$, we
have $\pi_{U^+_J}(ab)=\pi_{U^+_J}(\pi_{U^+}(a)b)$. Since $'U^+_J$
is a normal subgroup of $U^+$, $\pi_{U^+_J} \mid_{U^+}$ is a
homomorphism of $U^+$ onto $U^+_J$. Moreover, we have
$$\pi_{U^+_J}\bigl(x_i(a)\bigr)=\cases x_i(a), & \qquad \hbox{ if }
i \in J; \cr 1, & \qquad \hbox{ otherwise }. \cr
\endcases$$

Thus $\pi_{U^+_J}(U^+_{>0})=U^+_{w^J_0,
>0}$ and $\pi_{U^+_J}(U^+_{\geqslant 0})=U^+_{w^J_0,
\geqslant 0}$.

Let $U^-_J$  be the subgroup of $U^-$ generated by $\{U_{-\alpha}
\mid \alpha \in \Phi^+_J\}$ and $'U^-_J$ to be the subgroup of
$U^-$ generated by $\{U_{-\alpha} \mid \alpha \in
\Phi^+-\Phi^+_J\}$. Then we define $\pi_{U^-_J}: U^-@>>>U^-_J$ by
$\pi_{U^-_J}(u_1 u_2)=u_1$ for $u_1 \in U^-_J, u_2 \in 'U^-_J$.
(We will simply write $\pi_{U^-_I}$ as $\pi_{U^-}$.) We have
$\pi_{U^-_J}(U^-_{>0})=U^-_{w^J_0, >0}$ and $\pi_{U^-_J}(U^-_{\ge
0})=U^-_{w^J_0, \ge 0}$.

\subhead 1.8 \endsubhead For any vector space $V$ and a nonzero
element $v$ of $V$, we denote the image of $v$ in $P(V)$ by $[v]$.

If $(V, \rho)$ is a representation of $G$, we denote by $(V^*,
\rho^*)$ the dual representation of $G$. Then we have the standard
isomorphism $St_V: V \otimes V^*@>\simeq>>\End(V)$ defined by
$St_V(v \otimes v^*) (v')=v^* (v') v$ for all $v, v' \in V, v^*
\in V^*$. Now we have the $G \times G$ action on $V \otimes V^*$
by $(g_1, g_2) \cdot (v \otimes v^*)=(g_1 v) \otimes (g_2 v^*)$
for all $g_1, g_2 \in G, v \in V, v^* \in V^*$ and the $G \times
G$ action on $\End(V)$ by $\bigl((g_1, g_2) \cdot f \bigr)(v)=g_1
\bigl( f(g_2 \i v) \bigr)$ for all $g_1, g_2 \in G, f \in \End(V),
v \in V$. The standard isomorphism between $V \otimes V^*$ and
$\End(V)$ commutes with the $G \times G$ action. We will identify
$\End(V)$ with $V \otimes V^*$ via the standard isomorphism.

\head 2. The strata of the De Concini-Procesi Compactification
\endhead

\subhead 2.1 \endsubhead Let $\cv_G$ be the projective variety
whose points are the $\dim(G)$-dimensional Lie subalgebras of
$\Lie(G \times G)$. For any subset $J$ of $I$, define
$$Z_J=\{(P, Q, \gamma) \mid P \in \cp^J, Q \in \cp^{J^*},
\gamma=H_P g U_Q, P \opp ^g Q \}$$ with the $G \times G$ action by
$(g_1, g_2) \cdot (P, Q, H_P g U_Q)=\bigl( {^{g_1}P}, ^{g_2}Q,
H_{^{g_1}P} (g_1 g g_2\i) U_{^{g_2}Q} \bigr)$.

For $(P, Q, \gamma) \in Z_J$ and $g \in \gamma$, we set
$$H_{P, Q, \gamma}=\{(l+u_1, \Ad(g \i) l+u_2) \mid l \in \Lie(P
\cap {^g Q}), u_1 \in \Lie(U_P), u_2 \in \Lie(U_Q)\}.$$

Then $H_{P, Q, \gamma}$ is independent of the choice of $g$ (see
[L6, 12.2]) and is an element of $\cv_G$ (see [L6, 12.1]).
Moreover, $(P, Q, \gamma)@>>>H_{P, Q, \gamma}$ is an embedding of
$Z_J \subset \cv_G$ (see [L6, 12.2]). We will identify $Z_J$ with
the subvariety of $\cv_G$ defined above. Then we have
$\bar{G}=\bigsqcup\limits_{J \subset I} Z_J$, where $\bar{G}$ is
the the De Concini-Procesi compactification of $G$ (see [L6,
12.3]). We will call $\{Z_J \mid J \subset I\}$ the strata of
$\bar{G}$ and $Z_I$ (resp. $Z_{\varnothing}$) the highest (resp.
lowest) stratum of $\bar{G}$. It is easy to see that $Z_I$ is
isomorphic to $G$ and $Z_{\varnothing}$ is isomorphic to $\cb
\times \cb$.

Set $z^{\circ}_J=(P_J, Q_J, H_{P_J} U_{Q_J})$. Then $z^{\circ}_J
\in Z_J$ (see 1.4) and $Z_J=(G \times G) \cdot z^{\circ}_J$.

Since $G$ is adjoint, we have an isomorphism $\chi:
T@>\simeq>>(\bold R^*)^I$ defined by $\chi(t)=\bigl(\alpha_i(t) \i
\bigr)_{i \in I}$. We denote the closure of $T$ in $\bar{G}$ by
$\bar{T}$. We have $H_{P_J, Q_J, H_{P_J} U_{Q_J}}=\{(l+u_1, l+u_2)
\mid l \in \Lie(P_J \cap Q_J), u_1 \in U_{P_J}, u_2 \in U_{Q_J}
\}$. Moreover, for any $t \in Z(P_J \cap Q_J)$, $H_t$ is the
subspace of $\Lie(G) \times \Lie(G)$ spanned by the elements $(l,
l), (u_1, \Ad(t \i) u_1), (\Ad(t)u_2, u_2)$, where $l \in \Lie(P_J
\cap Q_J), u_1 \in U_{P_J}, u_2 \in U_{Q_J}$. Thus it is easy to
see that $z^{\circ}_J=\lim\limits_{\scriptstyle t_j=1,\forall j
\in J \atop \scriptstyle t_j@>>>0, \forall j \notin J} \chi \i
\bigl((t_i)_{i \in I} \bigr) \in \bar{T}$.

\proclaim{Proposition 2.2} The automorphism $\psi$ of the variety
$G$ (see 1.1) can be extended in a unique way to an automorphism
$\bar{\psi}$ of $\bar{G}$. Moreover, $\bar{\psi}(P, Q,
\gamma)=\bigl(\psi(Q), \psi(P), \psi(\gamma) \bigr) \in Z_J$ for
$J \subset I$ and $(P, Q, \gamma) \in Z_J$.
\endproclaim

Proof. The map $\psi: G@>>>G$ induces a bijective map $\psi:
\Lie(G)@>>>\Lie(G)$. Moreover, we have $\psi(\Ad(g) v)=\Ad
\bigl(\psi(g) \i \bigr) \psi(v)$ and $\psi(v+v')=\psi(v)+\psi(v')$
for $g \in G, v, v'\in \Lie(G)$. Now define $\delta: \Lie(G)
\times \Lie(G)@>>>\Lie(G) \times \Lie(G)$ by $\delta(v,
v')=\bigl(\psi(v'), \psi(v) \bigr)$ for $v, v' \in \Lie(G)$. Then
$\delta$ induces a bijection $\bar{\psi}: \cv_G@>>>\cv_G$.

Note that for any $g \in G$, we have $H_g=\{(v, \Ad(g) v) \mid v
\in \Lie{G}\}$ and $\bar{\psi}(H_g)=\{(\Ad(\psi(g) \i) \psi(v),
\psi(v)) \mid v \in \Lie(G)\}=H_{\psi(g)}$. Thus $\bar{\psi}$ is
an extension of the automorphism $\psi$ of $G$ into $\cv_G$.

Now for any $(P, Q, \gamma) \in Z_J$ and $g \in \gamma$, we have
$\psi(P) \in \cp^{J^*}, \psi(Q) \in \cp^J$ and $\psi(Q) \opp
^{\psi(g)} \psi(P)$ (see 1.4). Thus $\bigl(\psi(Q), \psi(P),
\psi(\gamma) \bigr) \in Z_J$. Moreover,
$$\eqalignno{\bar{\psi}(H_{P, Q, \gamma}) &=\{(\Ad(\psi(g))
\psi(l)+\psi(u_2), \psi(l)+\psi(u_1)) \mid l \in \Lie(P \cap {^g
Q}), \cr & \qquad \qquad \qquad \qquad \qquad \qquad \qquad \qquad
\qquad u_1 \in \Lie(U_P), u_2 \in \Lie(U_Q) \} \cr &=\{(l+u_2,
\Ad(\psi(g) \i)l+u_1) \mid l \in \Lie(\psi(Q) \cap ^{\psi(g)}
\psi(P)), \cr & \qquad \qquad \qquad \qquad \qquad \qquad \qquad
\qquad u_1 \in \Lie(\psi(U_P)), u_2 \in \Lie(\psi(U_Q)) \} \cr
&=H_{\psi(Q), \psi(P), \psi(\gamma)}. }$$

Thus $\bar{\psi} \mid_{\bar{G}}$ is an automorphism of $\bar{G}$.
Moreover, since $\bar{G}$ is the closure of $G$, $\bar{\psi}
\mid_{\bar{G}}$ is the unique automorphism of $\bar{G}$ that
extends the automorphism $\psi$ of $G$.

The proposition is proved. \qed

\subhead 2.3 \endsubhead For any $\lambda\in X$, set $\supp
(\lambda)=\{i \in I \mid<\alpha_i^\vee, \lambda> \neq 0 \}$.

In the rest of the section, I will fix a subset $J$ of $I$ and
$\lambda_1, \lambda_2 \in X^+$ with $\supp (\lambda_1)=I-J, \supp
(\lambda_2)=J$. Let $(V_{\lambda_1}, \rho_1)$ (resp.
$(V_{\lambda_2}, \rho_2)$) be the irreducible representation of
$G$ with the highest weight $\lambda_1$ (resp. $\lambda_2$).
Assume that $\dim V_{\lambda_1}=n_1, \dim V_{\lambda_2}=n_2$ and
$\{v_1, v_2, \ldots, v_{n_1}\}$ (resp. $\{v_1', v_2', \ldots,
v_{n_2}'\}$) is the canonical basis of $(V_{\lambda_1}, \rho_1)$
(resp. $(V_{\lambda_2}, \rho_2)$), where $v_1$ and $v_1'$ are the
highest weight vectors. Moreover, after reordering $\{2, 3,
\ldots, n_2\}$, we could assume that there exists some integer
$n_0 \in \{1, 2, \ldots, n_2\}$ such that for any $i \in \{1, 2,
\ldots, n_2\}$, the weight of $v_i'$ is of the form
$\lambda_2-\sum_{j \in J}a_j\alpha_j$ if and only if $i \le n_0$.

Define $i_J: G@>>>P \bigl(\End (V_{\lambda_1}) \bigr) \times P
\bigl(\End(V_{\lambda_2}) \bigr)$ by $i_J (g)=\Bigl( [\rho_1(g)],
[\rho_2(g)] \Bigr)$. Then since $\lambda_1+\lambda_2$ is a
dominant and regular weight, the closure of the image of $i_J$ in
$P \bigl(\End (V_{\lambda_1}) \bigr) \times P
\bigl(\End(V_{\lambda_2}) \bigr)$ is isomorphic to the De
Concini-Procesi compactification of $G$ (See [DP, 4.1]). We will
use $i_J$ as the embedding of $\bar{G}$ into $P \bigl(\End
(V_{\lambda_1}) \bigr) \times P \bigl(\End(V_{\lambda_2}) \bigr)$.
We will also identify $\bar{G}$ with its image under $i_J$.

\subhead 2.4 \endsubhead Now with respect to the canonical basis
of $V_{\lambda_1}$ and $V_{\lambda_2}$, we will identify
$\End(V_{\lambda_1})$ with $gl(n_1)$ and $\End(V_{\lambda_2})$
with $gl(n_2)$. Thus we will regard $\rho_1(g), \rho^*_1(g)$ as
$n_1 \times n_1$ matrices and $\rho_2(g), \rho^*_2(g)$ as $n_2
\times n_2$ matrices. It is easy to see that (in terms of
matrices) for any $g \in G, \rho^*_1(g)=^t \rho_1(g \i)$ and
$\rho^*_2(g)=^t \rho_2(g \i)$, where $^t M$ is the transpose of
the matrix $M$. Now for any $g_1, g_2 \in G$, $M_1 \in gl(n_1)$,
$M_2 \in gl(n_2)$, $(g_1, g_2) \cdot M_1=\rho_1(g_1) M_1
\rho_1(g_2 \i)$ and $(g_1, g_2) \cdot M_2=\rho_2(g_1) M_2
\rho_2(g_2 \i)$.

Set $L=P_J \cap Q_J$. Then $L$ is a reductive algebraic group with
the \'{e}pinglage ${(T, B^+ \cap L, B^- \cap L, x_j, y_j; j \in
J)}$. Now let $V_L$ be the subspace of $V_{\lambda_2}$ spanned by
$\{v_1', v_2', \ldots, v_{n_0}' \}$ and $I_L=(a_{ij}) \in
gl(n_2)$, where
$$a_{ij}=\cases 1, &\hbox{ if } i=j \in \{1, 2, \ldots, n_0\}; \cr
0, &\hbox{ otherwise}. \cr \endcases$$

Then $V_L$ is an irreducible representation of $L$ with the
highest weight $\lambda_2$ and canonical basis $\{v_1', v_2',
\ldots, v_{n_0}'\}$. Moreover, $\lambda_2$ is a dominant and
regular weight for $L$. Now set $I_1=\diag(1, 0, 0, \ldots, 0) \in
gl(n_1), I_2=\diag(1, 0, 0, \ldots, 0) \in gl(n_2)$. Then
$i_J(z^{\circ}_J)=\lim\limits_{\scriptstyle t_j=1,\forall j \in J
\atop \scriptstyle t_j@>>>0, \forall j \notin J} i_J \Bigl(\chi \i
\bigl((t_i)_{i \in I} \bigr) \Bigr)=\Bigl([v_1 \otimes v_1^*],
[\sum_{i=1}^{n_0} v'_i \otimes {v'_i}^*] \Bigr)=\Bigl( [I_1],
[I_L] \Bigr)$, where $\{{v_1}^*, {v_2}^*, \ldots, {v_{n_1}}^*\}$
(resp. $\{{v_1'}^*, {v_2'}^*, \ldots, {v_{n_2}'}^*\}$) is the dual
basis in $(V_{\lambda_1})^*$ (resp. $(V_{\lambda_2})^*$).

\subhead 2.5 \endsubhead Recall that $\supp(\lambda_1)=I-J$. Thus
for any $P \in \cp^J$, there is a unique $P$-stable line
$L_{\rho_1(P)}$ in $(V_{\lambda_1}, \rho_1)$ and $P \mapsto
L_{\rho_1(P)}$ is an embedding of $\cp^J$ into $P(V_{\lambda_1})$.
Similarly, for any $Q \in \cp^{J^*}$, there is a unique $Q$-stable
line $L_{\rho^*_1(Q)}$ in $(V^*_{\lambda_1}, \rho^*_1)$ and $Q
\mapsto L_{\rho^*_1(Q)}$ is an embedding of $\cp^{J^*}$ into
$P(V^*_{\lambda_1})$. It is easy to see $L_{\rho_1(P_J)}=[v_1]$,
$L_{\rho^*_1(Q_J)}=[{v_1}^*]$ and
$L_{\rho_1(^gP)}=\rho_1(g)L_{\rho_1(P)}, L_{\rho^*_1(^gQ)}
=\rho^*_1(g)L_{\rho^*_1(Q)}$ for $P \in \cp^J, Q \in \cp^{J^*}, g
\in G$.

There are projections $p_1: P \bigl(\End (V_{\lambda_1}) \bigr)
\times P \bigl(\End(V_{\lambda_2}) \bigr)@>>>P \bigl(\End
(V_{\lambda_1}) \bigr)$ and $p_2: P \bigl(\End (V_{\lambda_1})
\bigr) \times P \bigl(\End(V_{\lambda_2}) \bigr)@>>>P \bigl(\End
(V_{\lambda_2}) \bigr)$. It is easy to see that $p_1 \mid_{Z_J}$,
$p_2 \mid_{Z_J}$ commute with the $G \times G$ action and
$p_1(z^{\circ}_J)=[v_1 \otimes {v_1}^*]=[L_{\rho_1(P_J)} \otimes
L_{\rho^*_1(Q_J)}]$. Now for any $g_1, g_2 \in G$, we have
$$p_1\bigl((g_1, g_2) \cdot z^{\circ}_J\bigr)=[\rho_1(g_1)
L_{\rho_1(P_J)} \otimes \rho^*_1(g_2)
L_{\rho^*_1(Q_J)}]=[L_{\rho_1(^{g_1}P)} \otimes
L_{\rho^*_1(^{g_2}Q)}].$$ In other words, $p_1(z)=[L_{\rho_1(P)}
\otimes L_{\rho^*_1(Q)}]$ for $z=(P, Q, \gamma) \in Z_J$.

\subhead 2.6 \endsubhead Let $\overline{G_{>0}}$ be the closure of
$G_{>0}$ in $\bar G$. Then $\overline{G_{>0}}$ is also the closure
of $G_{\ge 0}$ in $\bar G$. We have $z^{\circ}_J \in
\overline{G_{>0}}$ (see 2.1). Now set
$$Z_{J, \ge 0}=Z_J \cap \overline{G_{>0}},$$
$$Z_{J, >0}=\{(g_1, g_2 \i) \cdot z^{\circ}_J \mid g_1, g_2 \in G_{>0}\}.$$

Since $\psi(G_{>0})=G_{>0}$, we have
$\bar{\psi}(\overline{G_{>0}})=\overline{G_{>0}}$. Moreover,
$\bar{\psi}(Z_J)=Z_J$ (see 2.2). Therefore $\bar{\psi}(Z_{J, \ge
0})=Z_{J, \ge 0}$. Similarly, $(g_1, g_2 \i) \cdot Z_{J, \ge 0}
\subset Z_{J, \ge 0}$ for any $g_1, g_2 \in G_{>0}$. Thus $Z_{J,
>0} \subset Z_{J, \ge 0}$. Moreover, it is easy to see that
$\bar{\psi}(Z_{J,>0})=Z_{J, >0}$.

Note that for any $u_1, u_4 \in U^-_{>0}, u_2, u_3 \in U^+_{>0},
t, t' \in T_{>0}$, we have $$\eqalignno{(u_1 u_2 t, u_3 \i u_4 \i
t') \cdot z^{\circ}_J&=(u_1 u_2, u_3 \i u_4 \i) \cdot (P_J, Q_J,
H_{P_J} t t' U_{Q_J})\cr &=(u_1, u_3 \i) \cdot \bigl(P_J, Q_J,
H_{P_J} \pi_{U^+_J}(u_2) t t' \pi_{U^-_J}(u_4) U_{Q_J} \bigr). }$$

Thus $$\eqalignno{Z_{J, >0} &=\{(u_1, u_2 \i) \cdot (P_J, Q_J,
H_{P_J} l U_{Q_J}) \mid u_1\in U^-_{>0}, u_2 \in U^+_{>0}, l \in
L_{>0}\} \cr &=\{(u'_1 t, {u'_2} \i) \cdot z^{\circ}_J \mid u'_1
\in U^-_{>0}, u'_2 \in U^+_{>0}, t \in T_{>0}\}. }$$

Moreover, for any $u_1, u'_1 \in U^-, u_2, u'_2 \in U^+$ and $t,
t' \in T$, it is easy to see that $(u_1 t, u_2) \cdot
z^{\circ}_J=(u'_1 t', u'_2) \cdot z^{\circ}_J$ if and only if
$(u_1 t) \i u'_1 t' \in l H_{P_J} \bigcap B^- \subset l Z(L)$ and
$u_2 \i u'_2 \in l \i H_{Q_J} \bigcap U^+ \subset l Z(L)$ for some
$l \in L$ , that is, $l \in Z(L)$, $u_1=u'_1, u_2=u'_2$ and $t \in
t' Z(L)$. Thus, $Z_{J,
>0} \cong U^-_{>0} \times U^+_{>0} \times T_{>0}/\bigl(T_{>0}
\bigcap Z(L) \bigr) \cong R_{>0}^{2 l(w_0)+|J|}$.

Now I will prove a criterion for $Z_{J, >0}$.

\proclaim{Theorem 2.7} Assume that $G$ is simply-laced. Let $z \in
Z_{J, \ge 0}$. Then $z \in Z_{J, >0}$ if and only if $z$ satisfies
the condition (*): $i_J(z)=\Bigl( [M_1], [M_2] \Bigr)$ and $i_J
\bigl(\bar{\psi}(z) \bigr)=\Bigl( [M_3], [M_4] \Bigr)$ for some
matrices $M_1, M_3 \in gl(n_1)$ and $M_2, M_4 \in gl(n_2)$ with
all the entries in $\bold R_{>0}$.
\endproclaim

Proof. If $z \in Z_{J, >0}$, then $z=(g_1, g_2 \i) \cdot
z^{\circ}_J$, for some $g_1, g_2 \in G_{>0}$. Assume that $g_1
\cdot v_1=\sum_{i=1}^{n_1} a_i v_i$ and $g_2 \i \cdot
v_1^*=\sum_{i=1}^{n_1} b_i v_i^*$. Then for any $i=1, 2, \ldots,
n_1$, $a_i, b_i>0$. Set $a_{ij}=a_i b_j$. Then
$p_1(z)=[\rho_1(g_1) I_1 \rho_1(g_2)]=[(a_{ij})]$ is a matrix with
all the entries in $\bold R_{>0}$.

We have $p_2(z)=[\rho_2(g_1) I_L \rho_2(g_2)]=[\rho_2(g_1) I_2
\rho_2(g_2)+ \rho_2(g_1) (I_L-I_2) \rho_2(g_2)]$. Note that
$\rho_2(g_1) I_2 \rho_2(g_2)$ is a matrix with all the entries in
$\bold R_{>0}$ and $\rho_2(g_1)$, $\rho_2(g_2)$, $(I_L-I_2)$ are
matrices with all the entries in $\bold R_{\ge 0}$. Thus
$\rho_2(g_1) (I_L-I_2) \rho_2(g_2)$ is a matrix with all its
entries in $\bold R_{\ge 0}$. So $\rho_2(g_1) I_L \rho_2(g_2)$ is
a matrix with all the entries in $\bold R_{>0}$.

Similarly, $i_J \bigl(\bar{\psi}(z) \bigr)=\Bigl( [M_3], [M_4]
\Bigr)$ for some matrices $M_3, M_4$ with all their entries in
$\bold R_{>0}$.

On the other hand, assume that $z$ satisfies the condition (*).
Suppose that $z=(P, Q, \gamma)$ and
$L_{\rho_1(P)}=[\sum_{i=1}^{n_1} a_i v_i]$,
$L_{\rho^*_1(Q)}=[\sum_{i=1}^{n_1} b_i v_i^*]$. We may also assume
that $a_{i_0}=b_{i_1}=1$ for some integers $i_0, i_1 \in \{1, 2,
\ldots, n_1\}$.

Set $M=(a_{ij}) \in gL(n_1)$, where $a_{ij}=a_i b_j$ for $i, j \in
\{1, 2, \ldots, n_1\}$. Then $p_1(z)=[L_{\rho_1(P)} \otimes
L_{\rho^*_1(Q)}]=[M]$. By the condition (*) and since $a_{i_0,
i_1}=a_{i_0} b_{i_1}=1$, we have that $M$ is a matrix with all its
entries in $\bold R_{>0}$. In particular, for any $i \in \{1, 2,
\ldots, n_1\}, a_{i, i_1}=a_i >0$. Therefore
$L_{\rho_1(P)}=[\sum_{i=1}^{n_1} a_i v_i]$, where $a_i>0$ for all
$i \in \{1, 2, \ldots, n_1\}$. By [R1, 5.1] (see also [L3, 3.4]),
$P \in \cp^J_{>0}$. Similarly, $\psi(Q) \in \cp^J_{>0}$. Thus
there exist $u_1 \in U^-_{>0}, u_2 \in U^+_{>0}$ and $l \in L$,
such that $z=(u_1, u_2 \i) \cdot (P_J, Q_J, H_{P_J} l U_{Q_J})$.

We can express $u_1, u_2$ in a unique way as $u_1=u_1'u_1''$, for
some $u_1' \in 'U^-_J$, $u_1'' \in U^-_J$ and $u_2=u_2''u_2'$, for
some $u_2' \in 'U^+_J$, $u_2'' \in U^+_J$ (see 1.7).

Recall that $V_L$ is the subspace of $V_{\lambda_2}$ spanned by
$\{v_1', v_2', \ldots, v_{n_0}' \}$. Let $V_L'$ be the subspace of
$V_{\lambda_2}$ spanned by $\{v_{n_0+1}', v_{n_0+2}', \ldots,
v_{n_2}' \}$. Then $u \cdot v-v \in V_L'$ and  $u \cdot V_L'
\subset V_L'$, for all $v \in V_L$, $\alpha \notin \Phi^+_J$ and
$u \in U_{-\alpha}$. Thus $u \cdot v-v \in V_L'$ and  $u \cdot
V_L' \subset V_L'$, for all $v \in V_L$ and $u \in 'U^-_J$.

Similarly, let $V_L^*$ be the subspace of $V^*_{\lambda_2}$
spanned by $\{{v'_1}^*, {v'_2}^*, \ldots, {v'_{n_0}}^*\}$ and
${V_L'}^*$ be the subspace of $V^*_{\lambda_2}$ spanned by
$\{{v'_{n_0+1}}^*, {v'_{n_0+2}}^*, \ldots, {v'_{n_2}}^*\}$. Then
for any  $v^* \in V_L^*$ and $u \in 'U^+_J$, we have $u \cdot v-v
\in {V_L'}^*$ and $u {V_L'}^* \subset {V_L'}^*$.

We define a map $\pi_L: gl(n_2)@>>>gl(n_0)$ by
$$\pi_L \bigl((a_{ij})_{i, j \in \{1, 2, \ldots, n_2\}} \bigr)=(a_{ij})_{i, j
\in \{1, 2, \ldots, n_0\}}$$

Then for any $u \in 'U^-_J, u' \in 'U^+_J$ and $M \in gl(n_2)$, we
have $\pi_L \bigl((u, u') \cdot M \bigr)=\pi_L(M)$. Set
$M_2=\rho_2(u_1 l) I_L \rho_2(u_2)$ and $l'=u_1'' l {u_2''} \in
L$. Then
$$\eqalignno{\pi_L(M_2)&=\pi_L \Bigl( (u_1, u_2 \i) \cdot \bigl(\rho_2(l)
I_L \bigr) \Bigr)=\pi_L \biggl( (u_1', {u_2'} \i) \cdot \Bigl(
(u_1'', {u_2''} \i) \cdot \bigl(\rho_2(l) I_L \bigr) \Bigr)
\biggr)\cr &=\pi_L \Bigl( (u_1'', {u_2''} \i) \cdot
\bigl(\rho_2(l) I_L \bigr) \Bigr)=\pi_L \bigl(\rho_2(l') I_L
\bigr)=\rho_L(l'). }$$

Since $p_2(z)=[M_2]$, $M_2$ is a matrix with all its entries
nonzero. Therefore $\rho_L(l')=\pi_L(M_2)$ is a matrix with all
its entries nonzero. Thus $l'=l_1 t_1 l_2$, for some $l_1 \in U^-
\cap L, l_2 \in U^+ \cap L, t_1 \in T$.

Set $\widetilde{u_1}=u_1' l_1$ and $\widetilde{u_2}=u_2' l_2$.
Then $^{\widetilde{u_1}}P_J=^{u_1 ({u_1''} \i l_1)} P_J
=^{u_1}P_J$. Similarly, we have $^{\widetilde{u_2} \i} Q_J=^{u_2
\i} Q_J$. So $z=(\widetilde{u_1}, \widetilde{u_2} \i) \cdot (P_J,
Q_J, H_{P_J} t_1 U_{Q_J})$.

Now for any $i_0, j_0 \in \{1, 2, \ldots, n_1\}$, define a map
$\pi^1_{i_0, j_0}: gl(n_1)@>>>\bold R$ by $\pi^1_{i_0, j_0}
\bigl((a_{ij})_{i, j \in \{1, 2, \ldots, n_1\}} \bigr)=a_{i_0,
j_0}$ and for any $i_0, j_0 \in \{1, 2, \ldots, n_2\}$, define a
map $\pi^2_{i_0, j_0}: gl(n_2)@>>>\bold R$ by $\pi^2_{i_0, j_0}
\bigl((a_{ij})_{i, j \in \{1, 2, \ldots, n_2\}} \bigr)=a_{i_0,
j_0}$.

Now $z=(\widetilde{u_1} t_1, \widetilde{u_2} \i) \cdot
z^{\circ}_J$ and $\bar{\psi}(z)=\bigl(\psi(\widetilde{u_2}) t_1,
\psi(\widetilde{u_1})\i \bigr) \cdot z^{\circ}_J$.

Set $$\leqalignno{\tilde{M_1}&=\rho_1(\widetilde{u_1} t_1) I_1
\rho_1(\widetilde{u_2}), \quad \tilde{M_3}=\rho_1
\bigl(\psi(\widetilde{u_2}) t_1 \bigr) I_1 \rho_1
\bigl(\psi(\widetilde{u_1}) \bigr), \cr
\tilde{M_2}&=\rho_2(\widetilde{u_1} t_1) I_L
\rho_2(\widetilde{u_2}), \quad \tilde{M_4}=\rho_2
\bigl(\psi(\widetilde{u_2}) t_1 \bigr) I_1 \rho_2
\bigl(\psi(\widetilde{u_1}) \bigr). \cr}$$

We have $\widetilde{u_1} \cdot v_1 =\sum_{i=1}^{n_1} {\pi^1_{i,
1}(\tilde{M_1}) \over \pi^1_{1, 1}(\tilde{M_1})} v_i$ and
$\psi(\widetilde{u_2}) \cdot v_1 =\sum_{i=1}^{n_1} {\pi^1_{i,
1}(\tilde{M_3}) \over \pi^1_{1, 1}(\tilde{M_3})} v_i$.

Moreover, let $V_0$ be the subspace of $V_{\lambda_2}$ spanned by
$\{v'_2, v'_3, \ldots, v'_{n_2}\}$ and ${V_0}^*$ be the subspace
of $V^*_{\lambda_2}$ spanned by $\{{v'_2}^*, {v'_3}^*, \ldots,
{v'_{n_2}}^*\}$. Then we have $u \cdot V_0 \subset V_0$, for all
$u \in U^-$ and $u' \cdot V^*_0 \subset V^*_0$, for all $u' \in
U^+$.

Thus for all $i=1, 2, \ldots, n_2$,
$$\eqalignno{\pi^2_{i, 1}(M_2) &=\pi^2_{i, 1} \bigl(\rho_2(\widetilde{u_1}
t_1) I_2 \rho_2(\widetilde{u_2}) \bigr)+\pi^2_{i, 1}
\bigl(\rho_2(\widetilde{u_1} t_1) (I_L-I_2)
\rho_2(\widetilde{u_2}) \bigr)\cr &=\pi^2_{i, 1}
\bigl(\rho_2(\widetilde{u_1} t_1) I_2 \rho_2(\widetilde{u_2})
\bigr).}$$

So $\widetilde{u_1} \cdot v_1' =\sum_{i=1}^{n_2} {\pi^2_{i,
1}(\tilde{M_2}) \over \pi^2_{1, 1}(\tilde{M_2})} v_i'$ and
$\psi(\widetilde{u_2}) \cdot v_1' =\sum_{i=1}^{n_2} {\pi^2_{i,
1}(\tilde{M_4}) \over \pi^2_{1, 1}(\tilde{M_4})} v_i'$. By [L2,
5.4], we have $\widetilde{u_1}, \psi(\widetilde{u_2}) \in
U^-_{>0}$. Therefore to prove that $z \in Z_{J, >0}$, it is enough
to prove that $t_1 \in T_{>0} Z(L)$, where $Z(L)$ is the center of
$L$.

For any $g \in (U^-, U^+) \cdot \bar{T}$, $g$ can be expressed in
a unique way as $g=(u_1, u_2) \cdot t$, for some $u_1 \in U^-$,
$u_2 \in U^+$, $t \in \bar{T}$. Now define $\pi_{\bar{T}}: (U^-,
U^+) \cdot \bar{T}@>>>\bar{T}$ by $\pi_{\bar{T}} \bigl( (u_1, u_2)
\cdot t \bigr)=t$ for all $u_1 \in U^-, u_2 \in U^+, t \in
\bar{T}$. Note that $(U^-, U^+) \cdot \bar{T} \cap
\overline{G_{>0}}$ is the closure of $G_{>0}$ in $(U^-, U^+) \cdot
\bar{T}$. Then $\pi_{\bar{T}}\bigl((U^-, U^+) \cdot \bar{T} \cap
\overline{G_{>0}} \bigr)$ is contained in the closure of $T_{>0}$
in $\bar{T}$. In particular, $\pi_{\bar{T}}(z)=t_1 t_J$ is
contained in the closure of $T_{>0}$ in $\bar{T}$. Therefore for
any $j \in J$, $\alpha_j(t_1)>0$. Now let $t_2$ be the unique
element in $T$ such that
$$\alpha_j(t_2)=\cases \alpha_j(t_1), &\hbox{ if } j \in J; \cr
\alpha_j(t_1)^2, &\hbox{ if } j \notin J. \cr \endcases$$

Then $t_2 \in T_{>0}$ and $t_2 \i t_1 \in Z(L)$. The theorem is
proved. \qed

\subhead Remark \endsubhead Theorem 2.7 is analogous to the
following statement in [L4, 5.4]: Assume that $G$ is simply laced
and $V$ is the irreducible representation of $G$ with the highest
weight $\lambda$, where $\lambda$ is a dominant and regular weight
of $G$. For any $g \in G$, let $M(g)$ be the matrix of $g: V@>>>V$
with respect to the canonical basis of $V$. Then for any $g \in
G$, $g \in G_{>0}$ if and only if $M(g)$ and $M \bigl(\psi(g)
\bigr)$ are matrices with all the entries in $\bold R_{>0}$.

\subhead 2.8 \endsubhead Before proving corollary 2.9, I will
introduce some technical tools.

Since $G$ is adjoint, there exists (in an essentially unique way)
$\tilde{G}$ with the \'{e}pinglage $(\tilde{T}, \tilde{B}^+,
\tilde{B}^-, \tilde{x}_{\tilde{i}}, \tilde{y}_{\tilde{i}};
\tilde{i} \in \tilde{I})$ and an automorphism $\sigma:
\tilde{G}@>>>\tilde{G}$ (over $\bold R$) such that the following
conditions are satisfied.

(a) $\tilde{G}$ is connected semisimple adjoint algebraic group
defined and split over $\bold R$.

(b) $\tilde{G}$ is simply laced.

(c) $\sigma$ preserves the \'{e}pinglage, that is,
$\sigma(\tilde{T})=\tilde{T}$ and there exists a permutation
$\tilde{i}@>>>\sigma(\tilde{i})$ of $\tilde{I}$, such that
$\sigma\bigl(\tilde{x}_{\tilde{i}}(a)
\bigr)=\tilde{x}_{\sigma(\tilde{i})}(a),
\sigma\bigl(\tilde{y}_{\tilde{i}}(a)
\bigr)=\tilde{y}_{\sigma(\tilde{i})}(a)$ for all $\tilde{i} \in
\tilde{I}$ and $a \in \bold R$.

(d) If $\tilde{i}_1 \neq \tilde{i}_2$ are in the same orbit of
$\sigma: \tilde{I}@>>>\tilde{I}$, then $\tilde{i}_1, \tilde{i}_2$
do not form an edge of the Coxeter graph.

(e) $\tilde{i}$ and $\sigma(\tilde{i})$ are in the same connected
component of the Coxeter graph, for any $\tilde{i} \in \tilde{I}$.

(f) There exists an isomorphism $\phi: \tilde{G}^{\sigma}@>>>G$
(as algebraic groups over $\bold R$) which is compatible with the
\'{e}pinglage of $G$ and the \'{e}pinglage $(\tilde{T}^{\sigma},
\tilde{B}^{+\sigma}, \tilde{B}^{-\sigma}, \tilde{x}_p,
\tilde{y}_p; p \in \bar{I})$ of $\tilde{G}^{\sigma}$, where
$\bar{I}$ is the set of orbit of $\sigma: \tilde{I}@>>>\tilde{I}$
and $\tilde{x}_p(a)=\prod\limits_{\tilde{i} \in p}
\tilde{x}_{\tilde{i}}(a), \tilde{y}_p(a)=\prod\limits_{\tilde{i}
\in p} \tilde{y}_{\tilde{i}}(a)$ for all $p \in \bar{I}$ and $a
\in \bold R$.

Let $\lambda$ be a dominant and regular weight of $\tilde{G}$ and
$(V, \rho)$ be the irreducible representation of $\tilde{G}$ with
highest weight $\lambda$. Let $\overline{\tilde{G}}$ be the
closure of $\{[\rho(\tilde{g})] \mid \tilde{g} \in \tilde{G}\}$ in
$P\bigl(\End(V)\bigr)$ and $\overline{\tilde{G}^{\sigma}}$ be the
closure of $\{[\rho(\tilde{g})] \mid \tilde{g} \in
\tilde{G}^{\sigma}\}$ in $P\bigl(\End(V)\bigr)$. Then since
$\lambda$ is a dominant and regular weight of $\tilde{G}$ and
$\lambda \mid_{\tilde{T}^{\sigma}}$ is a dominant and regular
weight of $\tilde{G}^{\sigma}$, we have that
$\overline{\tilde{G}}$ is the De Concini-Procesi compactification
of $\tilde{G}$ and $\overline{\tilde{G}^{\sigma}}$ is the De
Concini-Procesi compactification of $\tilde{G}^{\sigma}$. Since
$\overline{\tilde{G}}$ is closed in $P\bigl(\End(V)\bigr)$,
$\overline{\tilde{G}^{\sigma}}$ is the closure of
$\{[\rho(\tilde{g})] \mid \tilde{g} \in \tilde{G}^{\sigma}\}$ in
$\overline{\tilde{G}}$.

We have $\overline{\tilde{G}}=\bigsqcup\limits_{\tilde{J} \subset
\tilde{I}} \tilde{Z}_{\tilde{J}}=\bigsqcup\limits_{\tilde{J}
\subset \tilde{I}} (\tilde{G} \times \tilde{G}) \cdot
\tilde{z}^{\circ}_{\tilde{J}}$ and
$\overline{\tilde{G}^{\sigma}}=\bigsqcup\limits_{\tilde{J} \subset
\tilde{I}, \sigma{\tilde{J}}=\tilde{J}} (\tilde{G}^{\sigma} \times
\tilde{G}^{\sigma}) \cdot \tilde{z}^{\circ}_{\tilde{J}}$.
Moreover, $\sigma$ can be extended in a unique way to an
automorphism $\bar{\sigma}$ of $\overline{\tilde{G}}$. Since
$\overline{\tilde{G}}^{\bar{\sigma}}=\bigsqcup\limits_{\tilde{J}
\subset \tilde{I}, \sigma{\tilde{J}}=\tilde{J}}
(\tilde{Z}_{\tilde{J}})^{\bar{\sigma}}$ is a closed subset of
$\overline{\tilde{G}}$ containing $\tilde{G}^{\sigma}$, we have
$\overline{\tilde{G}^{\sigma}} \subset \bigsqcup\limits_{\tilde{J}
\subset \tilde{I}, \sigma{\tilde{J}}=\tilde{J}}
(\tilde{Z}_{\tilde{J}})^{\bar{\sigma}}$.

By the condition (f), there exists a bijection $\phi$ between
$\bar{I}$ and $I$, such that
$\phi\bigl(\tilde{x}_p(a)\bigr)=x_{\phi(p)}(a)$, for all $p \in
\bar{I}, a \in \bold R$. Moreover, the isomorphism $\phi$ from
$\tilde{G}^{\sigma}$ to $G$ can be extended in a unique way to an
isomorphism $\bar{\phi}:\overline{\tilde{G}^{\sigma}}@>>>\bar{G}$.
It is easy to see that for any $\tilde{J} \subset \tilde{I}$ with
$\sigma{\tilde{J}}=\tilde{J}$, we have
$\bar{\phi}\bigl((\tilde{G}^{\sigma} \times \tilde{G}^{\sigma})
\cdot \tilde{z}^{\circ}_{\tilde{J}}\bigr)=Z_{\phi \circ
\pi(\tilde{J})}$, where $\pi: \tilde{I}@>>>\bar{I}$ is the map
sending element of $\tilde{I}$ into the $\sigma$-orbit that
contains it.

\proclaim{Corollary 2.9} $Z_{J, \ge 0}=\bigcap\limits_{g_1, g_2
\in G_{>0}} (g_1 \i, g_2) \cdot Z_{J, >0}$ is the closure of
$Z_{J, >0}$ in $Z_J$. As a consequence, $Z_{J, \ge 0}$ and
$\overline{G_{>0}}$ are contractible.
\endproclaim

Proof. I will prove that $Z_{J, \ge 0} \subset
\bigcap\limits_{g_1, g_2 \in G_{>0}} (g_1 \i, g_2) \cdot Z_{J,
>0}$.

First, assume that $G$ is simply laced.

For any $g \in G_{>0}$, $i_J(g)=\Bigl( [\rho_1(g)], [\rho_2(g)]
\Bigr)$, where $\rho_1(g)$ and $\rho_2(g)$ are matrices with all
the entries in $\bold R_{>0}$. Then for any $z \in Z_{J, \ge 0}$,
we have $i_J(z)=\Bigl( [M_1], [M_2] \Bigr)$ for some matrices with
all the entries in $\bold R_{\ge 0}$. Similarly, $i_J
\bigl(\bar{\psi}(z) \bigr)=\Bigl( [M_3], [M_4] \Bigr)$ for some
matrices with all their entries in $\bold R_{\ge 0}$.

Note that for any $M_1', M_2', M_3' \in gl(n)$ such that $M_1',
M_3'$ are matrices with all their entries in $\bold R_{>0}$ and
$M_2'$ is a nonzero matrix with all the entries in $\bold R_{\ge
0}$, we have that $M_1' M_2' M_3'$ is a matrix with all the
entries in $\bold R_{>0}$. Thus for any $g_1, g_2 \in G_{>0}$, we
have that $(g_1, g_2 \i) \cdot z$ satisfies the condition (*) in
2.7. Moreover, $(g_1, g_2 \i) \cdot z \in Z_{J, \ge 0}$. Therefore
by 2.7, $(g_1, g_2 \i) \cdot z \in Z_{J, >0}$ for all $g_1, g_2
\in G_{>0}$.

In the general case, we will keep the notation of 2.8. Since the
isomorphism $\phi: \tilde{G}^{\sigma}@>>>G$ is compatible with the
\'{e}pinglages, we have
$\phi\bigl((\tilde{U}^{\pm}_{>0})^{\sigma}\bigr)=U^{\pm}_{>0}$,
$\phi\bigl((\tilde{T}_{>0})^{\sigma}\bigr)=T_{>0}$ and
$\phi\bigl((\tilde{G}_{>0})^{\sigma}\bigr)=G_{>0}$. Now for any $z
\in Z_{J, \ge 0}$, $z$ is contained in the closure of $G_{>0}$ in
$\bar{G}$. Thus $\bar{\phi}\i(z)$ is contained in the closure of
$(\tilde{G}_{>0})^{\sigma}$ in $\overline{\tilde{G}^{\sigma}}$,
hence contained in the closure of $(\tilde{G}_{>0})^{\sigma}$ in
$\overline{\tilde{G}}$. Therefore, $\bar{\phi}\i(z) \in
\tilde{Z}_{\tilde{J}, \ge0}$, where $\tilde{J}=\pi\i \circ
\phi\i(J)$.

For any $\widetilde{g_1}, \widetilde{g_2} \in
(\tilde{G}_{>0})^{\sigma}$, we have $(\widetilde{g_1},
\widetilde{g_2} \i) \cdot \bar{\phi}\i(z)=(\widetilde{u_1}
\tilde{t}, \widetilde{u_2}\i) \cdot \tilde{z}^{\circ}_{\tilde{J}}$
for some $\widetilde{u_1} \in \tilde{U}^-_{>0}, \widetilde{u_2}
\in \tilde{U}^+_{>0}, \tilde{t} \in \tilde{T}_{>0}$. Since
$\bar{\phi}\i(z) \in (\overline{\tilde{G}})^{\bar{\sigma}}$, we
have $(\widetilde{g_1}, \widetilde{g_2} \i) \cdot \bar{\phi}\i(z)
\in (\tilde{Z}_{\tilde{J}, >0})^{\bar{\sigma}}$. Then
$$\eqalignno{\bar{\sigma}\bigl((\widetilde{u_1} \tilde{t}, \widetilde{u_2}\i)
\cdot \tilde{z}^{\circ}_{\tilde{J}}\bigr)
&=\bigl(\sigma(\widetilde{u_1} \tilde{t}),
\sigma(\widetilde{u_2}\i)\bigr) \cdot
\bar{\sigma}(\tilde{z}^{\circ}_{\tilde{J}})=\bigl(\sigma(\widetilde{u_1})
\sigma(\tilde{t}), \sigma(\widetilde{u_2}\i)\bigr) \cdot
\tilde{z}^{\circ}_{\tilde{J}} \cr &=(\widetilde{u_1} \tilde{t},
\widetilde{u_2}\i) \cdot \tilde{z}^{\circ}_{\tilde{J}}.}$$

Thus $\sigma(\widetilde{u_1})=\widetilde{u_1}$ and
$\sigma(\widetilde{u_2})=\widetilde{u_2}$. Moreover, $(\tilde{t},
1) \cdot \tilde{z}^{\circ}_{\tilde{J}}=\bigl(\sigma(\tilde{t}), 1
\bigr) \cdot \tilde{z}^{\circ}_{\tilde{J}}$, that is,
$\tilde{\alpha}_{\tilde{j}}(\tilde{t})=\tilde{\alpha}_
{\tilde{j}}\bigl(\sigma((\tilde{t})\bigr)=\tilde{\alpha}_{\sigma(\tilde{j})}(\tilde{t})$
for all $\tilde{j} \in \tilde{J}$, where $\{\tilde{\alpha}_
{\tilde{i}} \mid \tilde{i} \in \tilde{I}\}$ is the set of simple
roots of $\tilde{G}$. Let $\tilde{t}'$ be the unique element in
$\tilde{T}$ such that
$$\tilde{\alpha}_ {\tilde{j}}(\tilde{t}')=\cases \tilde{\alpha}_
{\tilde{j}}(\tilde{t}), & \hbox{ if } \tilde{j} \in \tilde{J}; \cr
1, & \hbox{ otherwise }. \endcases$$

Then $\tilde{t}' \in (\tilde{T}_{>0})^{\sigma}$ and $(\tilde{t},
1) \cdot \tilde{z}^{\circ}_{\tilde{J}}=(\tilde{t}', 1) \cdot
\tilde{z}^{\circ}_{\tilde{J}}$. Thus $(\widetilde{g_1},
\widetilde{g_2} \i) \cdot \bar{\phi}\i(z)=(\widetilde{u_1}
\tilde{t}', \widetilde{u_2}\i) \cdot
\tilde{z}^{\circ}_{\tilde{J}}$. We have

$$\eqalignno{\bigl(\phi(\widetilde{g_1}), \phi(\widetilde{g_2})\i\bigr) \cdot
z &=\bar{\phi} \bigl((\widetilde{g_1}, \widetilde{g_2} \i) \cdot
\bar{\phi}\i(z)\bigr)=\bar{\phi} \bigl((\widetilde{u_1}
\tilde{t}', \widetilde{u_2}\i) \cdot
\tilde{z}^{\circ}_{\tilde{J}}\bigr)\cr
&=\bigl(\phi(\widetilde{u_1}) \phi(\tilde{t}'),
\phi(\widetilde{u_2}\i)\bigr) \cdot z^{\circ}_J \in Z_{J, >0}.}$$

Since $\phi\bigl((\tilde{G}_{>0})^{\sigma}\bigr)=G_{>0}$, we have
$Z_{J, \ge 0} \subset \bigcap\limits_{g_1, g_2 \in G_{>0}} (g_1
\i, g_2) \cdot Z_{J, >0}$.

Note that $(1, 1)$ is contained in the closure of $\{(g_1, g_2 \i)
\mid g_1, g_2 \in G_{>0} \}$. Hence, for any $z \in
\bigcap\limits_{g_1, g_2 \in G_{>0}} (g_1 \i, g_2) \cdot Z_{J,
>0}$, $z$ is contained in the closure of $Z_{J, >0}$. On the other
hand, $Z_{J, \ge 0}$ is a closed subset in $Z_J$. $Z_{J, \ge 0}$
contains $Z_{J,>0}$, hence contains the closure of $Z_{J, >0}$ in
$Z_J$. Therefore, $Z_{J, \ge 0}=\bigcap\limits_{g_1, g_2 \in
G_{>0}} (g_1 \i, g_2) \cdot Z_{J, >0}$ is the closure of $Z_{J,
>0}$ in $Z_J$.

Now set $g_r=\exp \bigl(r\sum\limits_{i \in I} (e_i+f_i) \bigr)$,
where $e_i$ and $f_i$ are the Chevalley generators related to our
\'{e}pinglage by $x_i(1)=\exp(e_i)$ and $y_i(1)=\exp(f_i)$. Then
$g_r \in G_{>0}$ for $r \in \bold R_{>0}$ (see [L1, 5.9]). Define
$f: R_{\ge 0} \times Z_{J, \ge 0}@>>>Z_{J, \ge 0}$ by $f(r,
z)=(g_r, g_r \i) \cdot z$ for $r \in R_{\ge 0}$ and $z \in Z_{J,
\ge 0}$. Then $f(0, z)=z$ and $f(1, z) \in Z_{J, >0}$ for all $z
\in Z_{J, \ge 0}$. Using the fact that $Z_{J, >0}$ is a cell (see
2.6), it follows that $Z_{J, \ge 0}$ is contractible.

Similarly, define $f': R_{\ge 0} \times
\overline{G_{>0}}@>>>\overline{G_{>0}}$ by $f'(r, z)=(g_r, g_r \i)
\cdot z$ for $r \in R_{\ge 0}$ and $z \in \overline{G_{>0}}$. Then
$f'(0, z)=z$ and $f'(1, z) \in \bigsqcup\limits_{K \subset I}
Z_{K, >0}$ for all $z \in \overline{G_{>0}}$. Note that
$\bigsqcup\limits_{K \subset I} Z_{K, >0}=\bigl(U^-_{>0},
(U^+_{>0})\i \bigr) \cdot \bigsqcup\limits_{K \subset I} (T_{>0},
1) \cdot z^{\circ}_K \cong U^-_{>0} \times U^+_{>0} \times
\bigsqcup\limits_{K \subset I} (T_{>0}, 1) \cdot z^{\circ}_K$ (see
2.6). Moreover, by [DP, 2.2], we have $\bigsqcup\limits_{K \subset
I} (T_{>0}, 1) \cdot z^{\circ}_K \cong R_{\ge 0}^I$. Thus
$\bigsqcup\limits_{K \subset I} Z_{K, >0} \cong R_{>0}^{2 l(w_0)}
\times R_{\ge 0}^I$ is contractible. Therefore $\overline{G_{>0}}$
is contractible. \qed

\head 3. The cell decomposition of $Z_{J, \ge 0}$ \endhead

\subhead 3.1 \endsubhead For any $P \in \cp^J, Q \in \cp^{J^*}, B
\in \cb$ and $g_1 \in H_P, g_2 \in U_Q, g \in G$, we have $\po
\bigl(P^B, ^{g_1 g g_2} (Q^B) \bigr)=\po \bigl({^{g_1 \i} (P^B)},
^{g g_2} (Q^B) \bigr)=\po \bigl(P^B, ^g (Q^B) \bigr)$. If
moreover, $P \opp ^g Q$, then $\po \bigl(P^B, ^g (Q^B) \bigr)=w
w_0$ for some $w \in W_J$ (see 1.4). Therefore, for any $v, v' \in
W$, $w, w' \in W^J$ and $y, y' \in W_J$ with $v \le w$ and $v' \le
w'$, Lusztig introduced the subset $Z^{v, w, v', w'; y, y'}_J$ and
$Z^{v, w, v', w'; y, y'}_{J, >0}$ of $Z_J$ which are defined as
follows.
$$\eqalignno{Z^{v, w, v', w'; y, y'}_J=\{(P, Q,& H_P g U_Q) \in Z_J \mid P \in
\cp^J_{v, w}, \psi(Q) \in \cp^J_{v', w'}, \cr &\po \bigl(P^{B^+},
^g (Q^{B^+}) \bigr)=y w_0, \po \bigl(P^{B^-}, ^g (Q^{B^-})
\bigr)=y' w_0\} }$$ and $$Z^{v, w, v', w'; y, y'}_{J, >0}=Z^{v, w,
v', w'; y, y'}_J \cap Z_{J, \ge 0}.$$

Then $$\eqalignno{Z_J &=\bigsqcup\limits_{\scriptstyle v, v' \in
W, w, w' \in W^J, y, y' \in W_J \atop \scriptstyle v \le w, v' \le
w'} Z^{v, w, v', w'; y, y'}_J, \cr Z_{J, \ge 0} &=
\bigsqcup\limits_{\scriptstyle v, v' \in W, w, w' \in W^J, y, y'
\in W_J \atop \scriptstyle v \le w, v' \le w'} Z^{v, w, v', w'; y,
y'}_{J, >0}. }$$

Lusztig conjectured that for any $v, v' \in W, w, w' \in W^J, y,
y' \in W_J$ such that $v \le w, v' \le w'$, $Z^{v, w, v', w'; y,
y'}_{J, >0}$ is either empty or a semi-algebraic cell. If it is
nonempty, then it is also a connected component of $Z^{v, w, v',
w'; y, y'}_J$.

In this section, we will prove this conjecture. Moreover, we will
show exactly when $Z^{v, w, v', w'; y, y'}_{J, >0}$ is nonempty
and we will give an explicit description of $Z^{v, w, v', w'; y,
y'}_{J, >0}$.

First, I will prove some elementary facts about the total
positivity of $G$.

\proclaim{Proposition 3.2} $$\eqalignno{\bigcap\limits_{u \in
U^{\pm}_{>0}} u \i U^{\pm}_{>0} &=\bigcap\limits_{u \in
U^{\pm}_{>0}} U^{\pm}_{>0} u \i=\bigcap\limits_{u \in U^{\pm}_{>
0}} u \i U^{\pm}_{\ge 0}=\bigcap\limits_{u \in U^{\pm}_{>0}}
U^{\pm}_{\ge 0} u \i=U^{\pm}_{\ge 0}, \cr \bigcap\limits_{g \in
G_{>0}} g \i G_{>0} &=\bigcap\limits_{g \in G_{>0}} G_{>0} g
\i=\bigcap\limits_{g \in G_{>0}} g \i G_{\ge 0}=\bigcap\limits_{g
\in G_{>0}} G_{\ge 0} g \i=G_{\ge 0}.}$$
\endproclaim

Proof. I will only prove $\bigcap\limits_{u \in U^+_{>0}} u \i
\cdot U^+_{>0}=U^+_{\ge 0}$. The rest of the equalities could be
proved in the same way.

Note that $u u_1 \in U^+_{>0}$ for all $u_1 \in U^+_{\ge 0}, u \in
U^+_{>0}$. Thus $u_1 \in \bigcap\limits_{u \in U^+_{>0}} u \i
\cdot U^+_{>0}$. On the other hand, assume that $u_1 \in
\bigcap\limits_{u \in U^+_{>0}} u \i \cdot U^+_{>0}$. Then $u u_1
\in U^+_{>0}$ for all $u \in U^+_{>0}$. We have
$u_1=\lim\limits_{\scriptstyle u \in U^+_{>0} \atop \scriptstyle
u@>>>1} u u_1$ is contained in the closure of $U^+_{>0}$ in $U^+$,
that is, $u_1 \in U^+_{\ge 0}$. So $\bigcap\limits_{u \in
U^+_{>0}} u \i \cdot U^+_{>0}=U^+_{\ge 0}$. \qed

For any $v, v' \in W$, $w, w' \in W^J$ such that $v \le w, v' \le
w'$, set $Z^{v, w, v', w'}_J =\bigsqcup\limits_{y, y' \in W_J}
Z^{v, w, v', w'; y, y'}_J$ and $Z^{v, w, v', w'}_{J,
>0}=\bigsqcup\limits_{y, y' \in W_J} Z^{v, w, v', w'; y, y'}_{J,
>0}$. We will give a characterization of $z \in Z^{v, w, v',
w'}_{J, >0}$ in 3.5.

\proclaim{Lemma 3.3} For any $w \in W$, $u \in U^-_{\ge 0}$, $\{
\pi_{U^+}(u_1 u) \mid u_1 \in U^+_{w, >0} \}=U^+_{w, >0}$.
\endproclaim

Proof. The following identities hold (see [L1, 1.3]):

(a) $t x_i(a)=x_i \bigl(\alpha_i(t) a \bigr)t, t y_i(a)=y_i
\bigl(\alpha_i(t) \i a \bigr)t$ for all $i \in I, t \in T, a \in
\bold R$.

(b) $y_{i_1}(a) x_{i_2}(b)=x_{i_2}(b) y_{i_1}(a)$ for all $a, b
\in \bold R$ and $i_1 \neq i_2 \in I$.

(c) $x_i(a) y_i(b)=y_i({b \over 1+ab}) \alpha_i^\vee({1 \over
1+ab}) x_i({a \over 1+ab})$ for all $a, b \in \bold R_{>0}, i \in
I$.

Thus $U^+_{w, >0} U^-_{\ge 0} \subset U^-_{\ge 0} T_{>0} U^+_{w,
>0}$ for $w \in W$. So we only need to prove that $U^+_{w, >0} \subset \{
\pi_{U^+}(u_1 u) \mid u_1 \in U^+_{w, >0} \}$. Now I will prove
the following statement:

$\{ \pi_{U^+} \bigl(u_1 y_i(a) \bigr) \mid u_1 \in U^+_{w, >0}
\}=U^+_{w, >0}$ for $i \in I, a \in \bold R_{>0}$.

We argue by induction on $l(w)$. It is easy to see that the
statement holds for $w=1$. Now assume that $w \neq 1$. Then there
exist $j \in I$ and $w_1 \in W$ such that $w=s_j w_1$ and
$l(w_1)=l(w)-1$. For any $u'_1 \in U^+_{w, >0}$, we have
$u'_1=u'_2 u'_3$ for some $u'_2 \in U^+_{s_j, >0}$ and $u'_3 \in
U^+_{w_1,
>0}$. By induction hypothesis, there exists $u_3 \in U^+_{w_1,
>0}, u' \in U^-$ and $t \in T$ such that $u_3 y_i(a)=u' t u'_3$.
Since $U^+_{w, >0} U^-_{s_i, >0} \subset U^-_{s_i, >0} T_{>0}
U^+_{w,
>0}$, we have $u' \in U^-_{s_i, >0}$ and $t \in T_{>0}$.

Now by (a), we have $t u'_2 t \i \in U^+_{s_j, >0}$. So by (b) and
(c), there exists $u_2 \in U^+_{s_j,
>0}$ such that $\pi_{U^+} (u_2 u')=t u'_2 t \i$. Thus
$$\eqalignno{\pi_{U^+} \bigl(u_2 u_3 y_i(a) \bigr)&=\pi_{U^+} \Bigl( (u_2 u')
\bigl( {u'}\i u_3 y_i(a) \bigr) \Bigr)=\pi_{U^+}
\bigl(\pi_{U^+}(u_2 u') {u'}\i u_3 y_i(a) \bigr)\cr &=\pi_{U^+}(t
u'_2 t \i t u'_3)=\pi_{U^+}(t u'_2 u'_3)=u'_1.}$$

So $u'_1 \in \{ \pi_{U^+}(u_1 y_i(a)) \mid u_1 \in U^+_{w, >0} \}
$. The statement is proved.

Now assume that $u \in U^-_{w', >0}$. I will prove the lemma by
induction on $l(w')$. It is easy to see that the lemma holds for
$w'=1$. Now assume that $w' \neq 1$. Then there exist $i \in I$
and $w'_1 \in W$ such that $l(w'_1)=l(w')-1$ and $w'=s_i w'_1$. We
have $u=y_i(a) u'$ for some $a \in \bold R_{>0}$ and $u' \in
U^-_{w'_1, >0}$. So
$$\eqalignno{\{ \pi_{U^+}(u_1 u) \mid u_1 \in U^+_{w, >0} \}&=
\{ \pi_{U^+} \bigl(u_1 y_i(a) u' \bigr) \mid u_1 \in U^+_{w, >0}
\}\cr &=\{ \pi_{U^+} \Bigl(\pi_{U^+} \bigl(u_1 y_i(a) \bigr) u
\Bigr) \mid u_1 \in U^+_{w, >0} \}\cr &=\{ \pi_{U^+}(u'_1 u') \mid
u'_1 \in U^+_{w, >0} \}. \cr}$$

By induction hypothesis, we have

$\{ \pi_{U^+}(u_1 u) \mid u_1 \in U^+_{w, >0} \}=\{ \pi_{U^+}(u'_1
u') \mid u'_1 \in U^+_{w, >0} \}=U^+_{w, >0}$. \qed

\proclaim{Lemma 3.4} Set $Z^1_{J, >0}=\{ (g_1, g_2 \i) \cdot
z^{\circ}_J \mid g_1 \in U^-_{\ge 0} T_{>0}, g_2 \in U^+_{\ge
0}\}$. Then

(a) $Z_{J, \ge 0}=\bigcap\limits_{u_1 \in U^+_{>0}, u_2 \i \in
U^-_{>0}} (u_1 \i, u_2) \cdot Z^1_{J, >0}$.

(b) $$\eqalignno{Z^1_{J, >0} &=\bigsqcup\limits_{w_1, w_2 \in W^J}
\{(^{u_1} P_J, ^{u_2 \i} Q_J, u_1 H_{P_J} l U_{Q_J} u_2) \mid u_1
\in U^-_{w_1, >0}, u_2 \in U^+_{w_2, >0}, l \in L_{\ge 0} \} \cr
&=\{ (P, Q, \gamma) \in Z_{J, \ge 0} \mid P=^{u_1} P_J,
\psi(Q)=^{u_2} P_J \hbox{ for some } u_1, u_2 \in U^-_{\ge 0}\}.
\cr}$$
\endproclaim

Proof. (a) By 2.9 and 3.2, we have

$$\eqalignno{Z_{J, \geqslant 0} &=\bigcap\limits_{g_1, g_2 \in
G_{>0}} (g_1 \i, g_2) \cdot Z_{J,
>0}=\bigcap\limits_{\scriptstyle t_1, t_2 \in T_{>0} \atop
\scriptstyle u_1, u_2 \in U^+_{>0} , u_3, u_4 \in U^-_{>0}} (u_1
\i u_3 \i t_1 \i, u_4 u_2 t_2) \cdot Z_{J,
>0}\cr &=\bigcap\limits_{u_1 \in U^+_{>0}, u_4 \in U^-_{>0}} (u_1
\i, u_4) \cdot \bigcap\limits_{u_2 \in U^+_{>0}, u_3 \in U^-_{>0}}
(u_2 \i, u_3) \cdot \bigcap\limits_{t_1, t_2 \in T_{>0}} (t_1 \i,
t_2) \cdot Z_{J,
>0}\cr &=\bigcap\limits_{u_1 \in U^+_{>0}, u_4 \in U^-_{>0}} (u_1
\i, u_4) \cdot \bigcap\limits_{u_2 \in U^+_{>0}, u_3 \in U^-_{>0}}
(u_2 \i, u_3) \cdot Z_{J, >0} \cr &=\bigcap\limits_{u_1 \in
U^+_{>0}, u_4 \in U^-_{>0}} (u_1 \i, u_4) \cdot
\bigcap\limits_{u_2 \in U^+_{>0}, u_3 \in U^-_{>0}} \bigl(u_2 \i
U^-_{>0} T_{>0}, (U^+_{>0} u_3\i)\i \bigr) \cdot z^{\circ}_J \cr
&=\bigcap\limits_{u_1 \in U^+_{>0}, u_2 \i \in U^-_{>0}} (u_1 \i,
u_2) \cdot \Bigl( \bigl(U^-_{\geqslant 0} T_{>0}, (U^+_{\geqslant
0}) \i \bigr) \cdot z^{\circ}_J \Bigr). \cr}$$

(b) For any $u \in U^-_{\ge 0}, v \in U^+_{\ge 0}, t \in T_{>0}$,
there exist $w_1, w_2 \in W^J, w_3, w_4 \in W_J$, such that $u=u_1
u_3$ for some $u_1 \in U^-_{w_1, >0}$, $u_3 \in U^-_{w_3, >0}$ and
$v=u_4 u_2$ for some $u_2 \in U^+_{w_2, >0}$, $u_4 \in U^+_{w_4,
>0}$. Then $(u t, v \i) \cdot
z^{\circ}_J=(^{u_1} P_J, ^{u_2 \i} Q_J, u_1 H_{P_J} u_3 t u_4
U_{Q_J} u_2)$. On the other hand, assume that $l \in L_{\ge 0}$,
then $l=u_3 t u_4$ for some $u_3 \in U^-_{\ge 0}, u_4 \in U^+_{\ge
0}, t \in T_{>0}$. Thus for any $u_1 \in U^-_{\ge 0}, u_2 \in
U^+_{\ge 0}$, we have

$(^{u_1} P_J, ^{u_2 \i} Q_J, u_1 H_{P_J} l U_{Q_J} u_2)=(u_1 u_3
t, u_2 \i u_4 \i) \cdot z^{\circ}_J \in Z^1_{J, >0}$.

Therefore $$\eqalignno{Z^1_{J, >0} &=\bigsqcup\limits_{w_1, w_2
\in W^J} \{(^{u_1} P_J, ^{u_2 \i} Q_J, u_1 H_{P_J} l U_{Q_J} u_2)
\mid u_1 \in U^-_{w_1, >0}, u_2 \in U^+_{w_2, >0}, l \in L_{\ge
0}\} \cr & \subset \{ (P, Q, \gamma) \in Z_{J, \ge 0} \mid
P=^{u_1} P_J, \psi(Q)=^{u_2} P_J \hbox{ for some } u_1, u_2 \in
U^-_{\ge 0}\}.}$$

Note that $\{^u P_J \mid u \in U^-_{\ge 0}\}=\bigsqcup\limits_{w
\in W^J} \{^u P_J \mid u \in U^-_{w, >0}\}$. Now assume that
$z=\bigl( {^{u_1} P_J}, ^{\psi(u_2) \i} Q_J, u_1 H_{P_J} l U_{Q_J}
\psi(u_2) \bigr)$ for some $w_1, w_2 \in W^J$ and $u_1 \in
U^-_{w_1, >0}$, $u_2 \in U^-_{w_2, >0}$, $l \in L$. To prove that
$z \in Z^1_{J, >0}$, it is enough to prove that $l \in L_{\ge 0}
Z(L)$. By part (a), for any $u_3, u_4 \in U^+_{>0}$,
$$\bigl(u_3, \psi(u_4) \i \bigr) \cdot z=\bigl( {^{u_3 u_1} P_J},
^{\psi(u_4 u_2) \i} Q_J, u_3 u_1 H_{P_J} l U_{Q_J} \psi(u_4 u_2)
\bigr) \in Z^1_{J, >0}.$$

Note that $u_3 u_1=u'_1 t_1 \pi_{U^+}(u_3 u_1)$ for some $u'_1 \in
U^-_{w_1, >0}, t_1 \in T_{>0}$ and $u_4 u_2=u'_2 t_2 \pi_{U^+}(u_4
u_2)$ for some $u'_2 \in U^-_{w_2, >0}, t_2 \in T_{>0}$. So we
have $^{u_3 u_1} P_J=^{u'_1} P_J$, $^{\psi(u_4 u_2) \i}
Q_J=^{\psi(u'_2) \i} Q_J$ and

$$\eqalignno{u_3 u_1 H_{P_J} l U_{Q_J} \psi(u_4 u_2) &=
u'_1 t_1 \pi_{U^+}(u_3 u_1) H_{P_J} l U_{Q_J} \psi
\bigl(\pi_{U^+}(u_4 u_2) \bigr) t_2 \psi(u'_2) \cr &=u'_1 H_{P_J}
t_1 \pi_{U^+_J}(u_3 u_1) l \psi\bigl(\pi_{U^+_J}(u_4 u_2) \bigr)
t_2 U_{Q_J} \psi(u'_2). \cr}$$

Then $t_1 \pi_{U^+_J}(u_3 u_1) l \psi\bigl(\pi_{U^+_J}(u_4 u_2)
\bigr) t_2 \in L_{\ge 0} Z(L)$. Since $t_1, t_2 \in T_{>0}$, we
have $\pi_{U^+_J}(u_3 u_1) l \psi\bigl(\pi_{U^+_J}(u_4 u_2) \bigr)
\in L_{\ge 0} Z(L)$ for all $u_3, u_4 \in U^+_{>0}$. By 1.8 and
3.3, $\pi_{U^+_J}(U^+_{>0} u_1)=\pi_{U^+_J} \bigl(
\pi_{U^+}(U^+_{>0} u_1) \bigr)=\pi_{U^+_J}(U^+_{>0})=U^+_{w^J_0,
>0}$. Similarly, we have $\pi_{U^+_J}(U^+_{>0} u_2)=U^+_{w^J_0, >0}$.
Thus $$\eqalignno{l & \in \bigcap\limits_{u_3, u_4 \in U^+_{w^J_0,
>0} } u_3 \i U^+_{w^J_0, \ge 0} T_{>0} Z(L) U^-_{w^J_0, \ge 0}
\psi(u_4) \i \cr &=U^+_{w^J_0, \ge 0} T_{>0} Z(L) U^-_{w^J_0, \ge
0}=L_{\ge 0} Z(L).}$$

The lemma is proved. \qed

\proclaim{Proposition 3.5} Let $z \in Z^{v, w, v', w'}_J$, then $z
\in Z^{v, w, v', w'}_{J, >0}$ if and only if for any $u_1 \in
U^+_{v \i, >0}, u_2 \in U^+_{{v'} \i, >0}, \bigl(u_1, \psi(u_2 \i)
\bigr) \cdot z \in Z^1_{J, >0}$.
\endproclaim

Proof. Assume that $z \in \bigcap\limits_{u_1 \in U^+_{v \i, >0},
u_2 \in U^+_{{v'} \i, >0}} \bigl(u_1 \i, \psi(u_2) \bigr) Z^1_{J,
>0}$. Then we have $z=\lim\limits_{u_1, u_2@>>>1} \bigl(u_1,
\psi(u_2) \i \bigr) \cdot z$ is contained in the closure of
$Z^1_{J, >0}$ in $Z_J$. Note that $Z_{J, >0} \subset Z^1_{J,
>0} \subset Z_{J, \ge 0}$. Thus by 2.9, $Z_{J, \ge 0}$ is the closure
of $Z^1_{J, >0}$ in $Z_J$. Therefore, $z$ is contained in $Z_{J,
\ge 0}$.

On the other hand, assume that $z=(P, Q, \gamma) \in Z^{v, w, v',
w'}_{J, >0}$. By 3.4(a), for any $u_1 \in U^+_{v \i, >0}$, $u_2
\in U^+_{{v'} \i, >0}$, we have $\bigl(u_1, \psi(u_2 \i) \bigr)
\cdot z \in Z_{J, \ge 0}$. Moreover, we have $^{u_1} P=^{u_1'}
P_J$ for some $u_1' \in U^-_{w, >0}$ (see 1.6). Similarly, we have
$\psi(^{\psi(u_2 \i)} Q)={^{u_2} \psi(Q)}=^{u_2'} P_J$ for some
$u_2' \in U^-_{w', >0}$. By 3.4(b), $\bigl(u_1, \psi(u_2 \i)
\bigr) \cdot z \in Z^1_{J, >0}$. \qed

\subhead 3.6 \endsubhead Now I will fix $w \in W^J$ and a reduced
expression $\hbox{\bf w}=(w_{(0)}, w_{(1)}, \ldots, w_{(n)})$ of
$w$. Assume that $w_{(j)}=w_{(j-1)} s_{i_j}$ for all $j=1, 2,
\ldots, n$. Let $v \le w$ and $\hbox{\bf v}_+=(v_{(0)},v_{(1)},
\ldots, v_{(n)})$ the positive subexpression of $\hbox{\bf w}$.

Define $$G_{\hbox{\bf v}_+, \hbox{\bf w}}=\Bigl\{ g=g_1 g_2 \cdots
g_k \Bigl | \hbox{$\matrix g_j=y_{i_j}(a_j) \text{ for } a_j \in
\bold R-\{0\}, & \text{ if } v_{(j-1)}=v_{(j)}\cr g_j
=\dot{s_{i_j}}, \hfill & \hbox{ if } v_{(j-1)}<v_{(j)}
\endmatrix$} \Bigr \},$$

$$G_{\hbox{\bf v}_+, \hbox{\bf w}, >0}=\Bigl\{ g=g_1 g_2
\cdots g_k \Bigl | \hbox{$\matrix g_j=y_{i_j}(a_j) \text{ for }
a_j \in \bold R_{>0}, & \text{ if } v_{(j-1)}=v_{(j)}\cr g_j
=\dot{s_{i_j}}, \hfill & \hbox{ if } v_{(j-1)}<v_{(j)}
\endmatrix$} \Bigr \}.$$

Marsh and Rietsch have proved that the morphism ${g \mapsto ^g
B^+}$ maps $G_{\hbox{\bf v}_+, \hbox{\bf w}}$ into $\car_{v, w}$
(see [MR, 5.2]) and $G_{\hbox{\bf v}_+, \hbox{\bf w}, >0}$
bijectively onto $\car_{v, w, >0}$ (see [MR, 11.3]).

The following proposition is a technical tool needed in the proof
of the main theorem.

\proclaim{Proposition 3.7} For any $g \in G_{\hbox{\bf v}_+,
\hbox{\bf w}, >0}$, we have
$$\bigcap\limits_{u \in U^+_{v \i,
>0}} \bigl(\pi_{U^+_J} (u g) \bigr) \i \cdot
U^+_{w^J_0, \ge 0}=\cases U^+_{w^J_0, \ge 0}, & \hbox{ if } v \in
W^J; \cr \varnothing, & \hbox{ otherwise}. \endcases$$
\endproclaim

The proof will be given in 3.13.

\proclaim{Lemma 3.8} Suppose $\alpha_{i_0}$ is a simple root such
that $v_1 \i \alpha_{i_0}>0$ for $v \le v_1 \le w$. Then for all
$g \in G_{\hbox{\bf v}_+, \hbox{\bf w}, >0}$ and $a \in \bold R$,
we have $x_{i_0}(a) g=g t g'$ for some $t \in T_{>0}$ and $g' \in
\prod\limits_{\a \in R(v)} U_{\a} \cdot \bigl({\dot{v}} \i
x_{i_0}(a) \dot{v} \bigr)$, where $R(v)=\{\a \in \Phi^+ \mid v \a
\in -\Phi^+\}$.
\endproclaim

Proof. Marsh and Rietsch proved in [MR, 11.8] that $g$ is of the
form $g=\bigl(\prod\limits_{j \in J^{\circ}_{\hbox{\bf v}_+}}
y_{v_{(j-1)} \alpha_{i_j}} (t_j) \bigr) \dot{v}$ and $v_{(j-1)}
\alpha_{i_1} \ne \alpha_{i_0}$, for all $j=1, 2, \ldots, n$. Thus
$g=g_1 \dot{v}$ for some $g_1 \in \prod\limits_{\alpha \in
\Phi^+-\{\alpha_{i_0}\}} U_{-\alpha}$. Set $T_1=\{t \in T \mid
\alpha_{i_0}(t)=1 \}$, then $T_1 \prod\limits_{\alpha \in
\Phi^+-\{\alpha_{i_0}\}} U_{-\alpha}$ is a normal subgroup of
$\psi(P_{\{i_0\}})$. Now set $x=x_{i_0}(a)$, then $x g_1 x \i \in
B^-$. We may assume that $x g_1 x \i=u_1 t_1$ for some $u_1 \in
U^-$ and $t_1 \in T$. Now $xg=x g_1 \dot{v}=(x g_1 x \i) x
\dot{v}=u_1 \dot{v} ({\dot{v}} \i t_1 \dot{v}) ({\dot{v}} \i x
\dot{v})$. Moreover, by [MR, 11.8], $x g \in g B^+$. Thus $x g=g_1
\dot{v} t_2 g_2 g_3=g_1 (\dot{v} t_2 g_2 t_2 \i \dot{v} \i)
\dot{v} t_2 g_3$, for some $t_2 \in T$, $g_2 \in \prod\limits_{\a
\in R(v)} U_{\alpha}$ and $g_3 \in \prod\limits_{\a \in
\Phi^+-R(v)} U_{\alpha}$. Note that $g_1 (\dot{v} t_2 g_2 t_2 \i
\dot{v} \i), u_1 \in U^-$, $t_2, {\dot{v}} \i t_1 \dot{v} \in T$
and $g_3, {\dot{v}} \i x \dot{v} \in \prod\limits_{\a \in
\Phi^+-R(v)} U_{\alpha}$. Thus $g_1 (\dot{v} t_2 g_2 t_2 \i
\dot{v} \i)=u_1$, $t_2={\dot{v}} \i t_1 \dot{v}$ and
$g_3={\dot{v}} \i x \dot{v}$. Note that $g \i x_{i_0}(b) g \in
B^+$ for $b \in \bold R$ (see [MR, 11.8]). We have that $\{\pi_T(g
\i x_{i_0}(b) g) \mid b \in \bold R\}$ is connected and contains
$\pi_T \bigl(g \i x_{i_0}(0) g \bigr)=1$. Hence $\pi_T(g \i
x_{i_0}(b) g) \in T_{>0}$ for $b \in \bold R$. In particular,
$\pi_T(g \i x g)=t_2 \in T_{>0}$. Therefore $x g=g t_2 g'$ with
$t_2 \in T_{>0}$ and $g'=g_2 g_3 \in \prod\limits_{\a \in R(v)}
U_{\a} \cdot ({\dot{v}} \i x \dot{v})$. \qed

\subhead Remark \endsubhead In [MR, 11.9], Marsh and Rietsch
pointed out that for any $j \in J^+_{\hbox{\bf v}_+}$, we have $u
\i \alpha_{i_j}>0$ for all $v_{(j)} \i v \le u \le w_{(j)} \i w$.

\subhead 3.9 \endsubhead Suppose that $J^+_{\hbox{\bf v}_+}=\{j_1,
j_2, \ldots, j_k\}$, where $j_1<j_2<\cdots<j_k$ and $g=g_1 g_2
\cdots g_n$, where $$g_j=\cases y_{i_j}(a_j) \hbox{ for } a_j \in
\bold R_{>0}, & \hbox{ if } j \in J^{\circ}_{\hbox{\bf v}_+}; \cr
\dot{s_{i_j}}, & \hbox{ if } j \in J^+_{\hbox{\bf v}_+}. \cr
\endcases$$

For any $m=1, \ldots, k$, define $v_m=v_{(j_m)} \i v$,
$g_{(m)}=g_{j_m+1} g_{j_m+2} \cdots g_n$ and $f_m(a)=g_{(m)} \i
x_{i_{j_m}}(-a) g_{(m)} \in B^+$ (see [MR, 11.8]). Now I will
prove the following lemma.

\proclaim{Lemma 3.10} Keep the notation in 3.9. Then

(a) For any $u \in U^+_{v \i, >0}$, $u g=g' t u'$ for some $g' \in
U^-_{w,
>0}, t \in T_{>0}$ and $u' \in U^+$.

(b) $$\pi_{U^+}(U^+_{v \i, >0} g)=\{\pi_{U^+} \bigl(f_k(a_k)
f_{k-1}(a_{k-1} \bigr) \cdots f_1(a_1) \bigr) \mid a_1, a_2,
\ldots, a_k \in \bold R_{>0}\}.$$

\endproclaim

Proof. I will prove the lemma by induction on $l(v)$. It is easy
to see that the lemma holds when $v=1$. Now assume that $v \ne 1$.

For any $u \in U^+_{v \i, >0}$, since $^g B^+ \in \car_{v, w,
>0}$, we have $^{u g} B^+ \in \car_{1, w, >0}$. Thus $u g=g' t u'$
for some $g' \in U^-_{w, >0}, t \in T$ and $u' \in U^+$. Set
$y=g_{i_1} g_{i_2} \cdots g_{i_{j_1-1}}$. Note that $y \in
U^-_{\ge 0}$, we have $u y=y' t u'$ for some $y' \in U^-$, $u' \in
U^+_{v \i, >0}$ and $t \in T_{>0}$. Hence $\pi_T(u g)=\pi_T(u y
\dot{s_{i_{j_1}}} g_{(1)})=\pi_T(y' t u' \dot{s_{i_{j_1}}}
g_{(1)}) \in T_{>0} \pi_T(u' \dot{s_{i_{j_1}}} g_{(1)})$. To prove
that $\pi_T(U^+_{v \i, >0} g) \subset T_{>0}$, it is enough to
prove that $\pi_T(u \dot{s_{i_{j_1}}} g_{(1)}) \in T_{>0}$ for all
$u \in U^+_{v \i, >0}$.

For any $u \in U^+_{v \i, >0}$, we have $u=u_1 x_{i_{j_1}}(a)$ for
some $u_1 \in U^+_{v \i s_{i_{j_1}}, >0}$ and $a \in \bold
R_{>0}$. It is easy to see that $x_{i_{j_1}}(a) \dot{s_{i_{j_1}}}
g_{(1)}=\alpha_{i_{j_1}}^\vee(a) y_{i_{j_1}}(a) x_{i_{j_1}}(-a \i)
g_{(1)}$. Note that $\alpha_{i_{j_1}}^\vee(a) \in T_{>0}$ and by
3.8, ${g_{(1)}}\i x_{i_{j_1}}(-a \i) g_{(1)} \in T_{>0} U^+$.
Hence by 1.7, we have
$$\eqalignno{\pi_T(u \dot{s_{i_{j_1}}} g_{(1)}) &=\pi_T \Bigl(u_1
\alpha_{i_{j_1}}^\vee(a) y_{i_{j_1}}(a) g_{(1)} \bigl({g_{(1)}}\i
x_{i_{j_1}}(-a \i) g_{(1)} \bigr) \Bigr) \cr & \in  T_{>0} \pi_T
\bigl(U^+_{v \i s_{i_{j_1}}, >0} y_{i_{j_1}}(a) g_{(1)} \bigr)
T_{>0}. \cr}$$

Set $$\eqalignno{\hbox{\bf w}' &=(1, w_{(j_1-1)} \i w_{(j_1)},
\ldots, w_{(j_1-1)} \i w_{(n)}), \cr \hbox{\bf v}'_+ &=(1,
s_{i_{j_1}} v_{(j_1)}, s_{i_{j_1}} v_{(j_1+1)}, \ldots,
s_{i_{j_1}} v_{(n)}). \cr}$$

Then $\hbox{\bf w}'$ is a reduced expression of $w_{(j_1-1)} \i
w_{(n)}$ and $\hbox{\bf v}'_+$ is a positive subexpression of
$\hbox{\bf w}'$. For any $a \in \bold R_{>0}$, $y_{i_{j_1}}(a)
g_{(1)} \in G_{\hbox{\bf v}'_+, \hbox{\bf w}', >0}$. Thus by
induction hypothesis, for any $a \in \bold R_{>0}$, $\pi_T(U^+_{v
\i s_{i_{j_1}}, >0} y_{i_{j_1}}(a) g_{(1)}) \subset T_{>0}$.
Therefore, $\pi_T(u g) \in T_{>0}$. Part (a) is proved.

We have $$\eqalignno{\pi_{U^+}(U^+_{v \i, >0} g)&=\pi_{U^+}(U^+_{v
\i, >0} y \dot{s_{i_{j_1}}} g_{(1)})=\pi_{U^+}(\pi_{U^+}
\bigl(U^+_{v \i, >0} y) \dot{s_{i_{j_1}}} g_{(1)} \bigr)\cr
&=\pi_{U^+}(U^+_{v \i, >0} \dot{s_{i_{j_1}}}
g_{(1)})=\bigcup\limits_{a \in \bold R_{>0}} \pi_{U^+}
\bigl(U^+_{v\i s_{i_{j_1}}, >0} x_{i_{j_1}}(a\i) \dot{s_{i_{j_1}}}
g_{(1)} \bigr) \cr &=\bigcup\limits_{a \in \bold R_{>0}} \pi_{U^+}
\bigl(U^+_{v \i s_{i_{j_1}}, >0} \alpha_{i_{j_1}}^\vee(a \i)
y_{i_{j_1}}(a \i) g_{(1)} f_1(a) \bigr) \cr &=\bigcup\limits_{a
\in \bold R_{>0}} \pi_{U^+} \Bigl(\pi_{U^+} \bigl(U^+_{v \i
s_{i_{j_1}}, >0} \alpha_{i_{j_1}}^\vee(a \i) y_{i_{j_1}}(a \i)
\bigr) g_{(1)} f_1(a) \Bigr) \cr &=\bigcup\limits_{a \in \bold
R_{>0}} \pi_{U^+} \bigl(U^+_{v \i s_{i_{j_1}}, >0} g_{(1)} f_1(a)
\bigr)=\bigcup\limits_{a \in \bold R_{>0}} \pi_{U^+}
\Bigl(\pi_{U^+}\bigl(U^+_{v \i s_{i_{j_1}}, >0} g_{(1)} \bigr)
f_1(a) \Bigr). \cr}$$

By induction hypothesis, $$\eqalignno{\pi_{U^+}(U^+_{v \i
s_{i_{j_1}}, >0} g_{(1)})=\{\pi_{U^+} \bigl(f_k(a_k)
f_{k-1}(a_{k-1}) \cdots f_2(a_2) \bigr) \mid a_2, a_3, \ldots, a_k
\in \bold R_{>0}\}.}$$

Thus $$\eqalignno{\pi_{U^+}(U^+_{v \i, >0} g) &=\bigcup\limits_{a
\in \bold R_{>0}} \pi_{U^+} \Bigl(\pi_{U^+}\bigl(U^+_{v \i
s_{i_{j_1}}, >0} g_{(1)} \bigr) f_1(a) \Bigr) \cr &=\{\pi_{U^+}
\bigl(f_k(a_k) f_{k-1}(a_{k-1}) \cdots f_1(a_1) \bigr) \mid a_1,
a_2, \ldots, a_k \in \bold R_{>0}\}. \qed}$$

\subhead Remark. \endsubhead The referee pointed out to me that
the assertion $t \in T_{>0}$ of 3.10(a) could also been proved
using generalized minors.

\proclaim{Lemma 3.11} Assume that $\alpha$ is a positive root and
$u \in U_{\alpha}$, $u' \in U^+$ such that $u^n u' \in U^+_{\ge
0}$ for all $n \in \bold N$. Then $u=x_i(a)$ for some $i \in I$
and $a \in \bold R_{\ge 0}$.
\endproclaim

Proof. There exists $t \in T_{>0}$, such that $\alpha_i(t)=2$ for
all $i \in I$. Then $t u t\i=u^{\alpha(t)}=u^m$ for some $m \in
\bold N$. By assumption, $t^n u t^{-n} u' \in U^+_{\ge 0}$ for all
$n \in \bold N$. Thus $u \bigl(t^{-n} u' t^n \bigr)=t^{-n}
\bigl(t^n u t^{-n} u' \bigr) t^n \in U^+_{\ge 0}$. Moreover, it is
easy to see that $\lim\limits_{n \rightarrow \infty} t^{-n} u'
t^n=1$. Since $U^+_{\ge 0}$ is a closed subset of $U^+$,
$\lim\limits_{n \rightarrow \infty} u t^{-n} u' t^n=u \in U^+_{\ge
0}$. Thus $u=x_i(a)$ for some $i \in I$ and $a \in \bold R_{\ge
0}$. \qed

\proclaim{Lemma 3.12} Assume that $w \in W$ and $i, j \in I$ such
that $w \i \alpha_i=\alpha_j$. Then there exists $c \in \bold
R_{>0}$, such that $\dot{w} \i x_i(a) \dot{w}=x_j(c a)$ for all $a
\in \bold R$.
\endproclaim

Proof. There exist $c, c' \in \bold R-\{0\}$, such that $y_i(a)
\dot{w}=\dot{w} y_j(c' a)$ and $x_i(a) \dot{w}=\dot{w} x_j(c a)$
for $a \in \bold R$. Since $^{\dot{w}} B^- \in \cb_{\ge 0}$, we
have $^{y_i(1) \dot{w}} B^+=^{\dot{w} y_j(c')} B^+ \in \cb_{\ge
0}$. By 3.6, $c' \ge 0$. Thus $c'>0$. Moreover, since $w
\alpha_j=\alpha_i>0$, we have $w s_j w \i=s_i$ and $l(w s_j)=l(s_i
w)=l(w)+1$. Hence, setting $w'=w s_j=s_i w$, we have
$\dot{w}'=\dot{w} \dot{s_j}=\dot{s_i} \dot{w}$, that is $\dot{w}
x_i(-1) y_i(1) x_i(-1)=x_j(-c) y_j(c') x_i(-c) \dot{w}=x_j(-1)
y_j(1) x_j(-1) \dot{w}$. Therefore, $c={c'} \i>0$. \qed

\subhead 3.13. Proof of Proposition 3.7 \endsubhead If $v \in
W^J$, then $v \alpha>0$ for $\alpha \in \Phi^+_J$. So
$\pi_{U^+_J}(\prod\limits_{\a \in R(v)} U_{\alpha})=\{1\}$. By
3.8, $f_m(a) \in T (\prod\limits_{\a \in R(v_m)} U_{\a}) \cdot
U_{v_m \i \a_{i_{j_m}}}$ for all $m \in \{1, 2, \ldots, k\}$. Note
that $v \a \in -\Phi^+$ for all $a \in R(v_m)$ and $v v_m \i
\a_{i_{j_m}}=v_{(j_m)} \a_{i_{j_m}} \in -\Phi^+$. So $f_m(a) \in T
\prod\limits_{\a \in R(v)} U_{\alpha}$ and $f_k(a_k)
f_{k-1}(a_{k-1}) \cdots f_1(a_1) \in T \prod\limits_{\a \in R(v)}
U_{\alpha}$. Hence by 3.10(b), $\pi_{U^+_J}(u g)=1$ for all $u \in
U^+_{v \i, >0}$. Therefore $\bigcap\limits_{u \in U^+_{v \i,
>0}} \bigl(\pi_{U^+_J} (u g) \bigr) \i \cdot U^+_{w^J_0, \ge
0}=U^+_{w^J_0, \ge 0}$.

If $v \notin W^J$, then there exists $\alpha \in \Phi^+_J$ such
that $v \a \in -\Phi^+_J$, that is, $v_m \i \alpha_{i_{j_m}} \in
\Phi^+_J$ for some $m \in \{1, 2, \ldots, k\}$. Set $k_0=\max\{m
\mid v_m \i \alpha_{i_{j_m}} \in \Phi^+_J\}$. Then since
$R(v_{k_0})=\{v_m \i \alpha_{i_{j_m}} \mid m>k_0\}$, we have that
$v_{k_0} \a>0$ for $\alpha \in \Phi^+_J$. Hence by 3.8,
$\pi_{U^+_J} \bigl(f_{k_0}(a) \bigr)=\dot{v_{k_0}} \i
x_{i_{j_{k_0}}}(-a) \dot{v_{k_0}}$. If $u' \in \bigcap\limits_{u
\in U^+_{v \i, >0}} \bigl(\pi_{U^+_J} (u g) \bigr) \i \cdot
U^+_{w^J_0, \ge 0}$, then $\pi_{U^+_J} \bigl(f_k(a_k)
f_{k-1}(a_{k-1}) \cdots f_1(a_1) \bigr) u' \in U^+_{w^J_0, \ge 0}$
for all $a_1, a_2, \ldots, a_k \in \bold R_{>0}$. Since
$U^+_{w^J_0, \ge 0}$ is a closed subset of $G$, $\pi_{U^+_J}
\bigl(f_k(a_k) f_{k-1}(a_{k-1}) \cdots f_1(a_1) \bigr) u' \in
U^+_{w^J_0, \ge 0}$ for all $a_1, a_2, \ldots, a_k \in \bold
R_{\ge 0}$. Now take $a_m=0$ for $m \in \{1, 2 , \ldots,
k\}-\{k_0\}$, then $\pi_{U^+_J} \bigl(f_{k_0}(a) \bigr) u' \in
U^+_{w^J_0, \ge 0}$ for all $a \in \bold R_{>0}$. Set
$u_1=\dot{v_{k_0}} \i x_{i_{j_{k_0}}}(-1) \dot{v_{k_0}}$. Then
$u_1^n u' \in U^+_{w^J_0, \ge 0}$ for all $n \in N$. Thus by 3.11,
$v_{k_0} \i \alpha_{i_{j_{k_0}}}=\alpha_{j'}$ for some $j' \in J$
and $u_1 \in U^+_{w^J_0, \ge 0}$. By 3.12, $u_1=x_{j'}(-c)$ for
some $c \in \bold R_{>0}$. That is a contradiction. The
proposition is proved. \qed

\

 Let me recall that $L=P_J \bigcap Q_J$ (see 2.4). Now I will prove the main theorem.

\proclaim{Theorem 3.14} For any $v, w, v', w' \in W^J$ such that
$v \le w, v' \le w'$, set $$\eqalignno{\tilde{Z}^{v, w, v',
w'}_{J, >0}=\Bigl \{\bigl(^g P_J, ^{\psi(g') \i}Q_J, g H_{p_J} l
U_{Q_J} \psi(g') \bigr) \Bigl | \hbox{$\matrix g \in G_{\hbox{\bf
v}_+, \hbox{\bf w}, >0}, & g' \in G_{\hbox{\bf v}'_+, \hbox{\bf
w}', >0} \cr \text{and } l \in L_{\ge 0} \hfill & \endmatrix $}
\Bigr\}.}$$

Then
$$Z^{v, w, v', w'}_{J, >0}=\cases \tilde{Z}^{v, w, v', w'}_{J,
>0}, & \hbox{ if } v, w, v', w' \in W^J, v \le w, v' \le w'; \cr
\varnothing, & \hbox{ otherwise. } \endcases$$
\endproclaim

Proof. Note that $\{(P, Q, \gamma) \in Z_J \mid P \in \cp^J_{\ge
0}, \psi(Q) \in \cp^J_{\ge 0}\}$ is a closed subset containing
$Z_{J, >0}$. Hence it contains $Z_{J, \ge 0}$. Now fix $g \in
G_{\hbox{\bf v}_+, \hbox{\bf w}, >0}, g' \in G_{\hbox{\bf v}'_+,
\hbox{\bf w}', >0}$ and $l \in L$. By 3.10 (a), for any $u \in
U^+_{v \i, >0}$, $u g=a t \pi_{U^+}(u g)$ for some $a \in U^-_{w,
>0}$ and $t \in T_{>0}$. Similarly, for any $u' \in U^+_{{v'} \i,
>0}$, $u' g'=a' t' \pi_{U^+}(u' g')$ for some $a' \in U^-_{w',
>0}$ and $t' \in T_{>0}$. Set $z=\bigl(^g P_J, ^{\psi(g') \i}Q_J,
g H_{p_J} l U_{Q_J} \psi(g') \bigr)$. We have
$$\eqalignno{ \bigl(u, \psi(u') \i \bigr) \cdot z &=
\Bigl({^a P_J}, ^{\psi(a') \i} Q_J, a t \pi_{U^+}(u g) H_{P_J} l
U_{Q_J} \psi \bigl(\pi_{U^+}(u' g') \bigr) t' \psi(a') \Bigr) \cr
&=\Bigl({^a P_J}, ^{\psi(a') \i} Q_J, a H_{P_J} t \pi_{U^+_J}(u g)
l \psi \bigl(\pi_{U^+_J}(u' g') \bigr) t' U_{Q_J} \psi(a') \Bigr).
\cr}$$

Then $\bigl(u, \psi(u') \i \bigr) \cdot z \in Z^1_{J, >0}$ if and
only if $t \pi_{U^+_J}(u g) l \psi \bigl(\pi_{U^+_J}(u' g') \bigr)
t' \in L_{\ge 0} Z(L)$, that is, $$\eqalignno{l & \in
\pi_{U^+_J}(u g) \i L_{\ge 0} Z(L) \psi \bigl(\pi_{U^+_J}(u' g')
\bigr) \i \cr &=\bigl(\pi_{U^+_J}(u g) \i U^+_{w^J_0, \ge 0}
\bigr) T_{>0} Z(L) \psi \bigl(\pi_{U^+_J}(u' g') \i U^+_{w^J_0,
\ge 0} \bigr). \cr}$$

So by 3.5, $z \in Z_{J, \ge 0}$ if and only if
$$\eqalignno{l & \in \bigcap\limits_{\scriptstyle u \in U^+_{v \i,
>0} \atop \scriptstyle u' \in U^+_{{v'} \i, >0}}
\bigl(\pi_{U^+_J}(u g) \i U^+_{w^J_0, \ge 0} \bigr) T_{>0} Z(L)
\psi \bigl(\pi_{U^+_J}(u' g') \i U^+_{w^J_0, \ge 0} \bigr) \cr
&=\bigcap\limits_{u \in U^+_{v \i, >0}} \bigl(\pi_{U^+_J}(u g) \i
U^+_{w^J_0, \ge 0} \bigr) T_{>0} Z(L) \psi
\bigl(\bigcap\limits_{u' \in U^+_{{v'} \i,
>0}} \pi_{U^+_J}(u' g') \i U^+_{w^J_0, \ge 0} \bigr). \cr}$$

By 3.7, $z \in Z_{J, \ge 0}$ if and only if $v, v' \in W^J$ and $l
\in L_{\ge 0} Z(L)$. The theorem is proved. \qed

\subhead 3.15 \endsubhead It is known that $G_{\ge 0}
=\bigsqcup\limits_{w, w' \in W} U^-_{w, >0} T_{>0} U^+_{w', >0}$,
where for any $w, w' \in W$, $U^-_{w, >0} T_{>0} U^+_{w', >0}$ is
a semi-algebraic cell (see [L1, 2.11]) and is a connected
component of $B^+ \dot{w} B^+ \cap B^- \dot{w}' B^-$ (see [FZ]).
Moreover, Rietsch proved in [R2, 2.8] that $\cb_{\ge 0}
=\bigsqcup\limits_{v \le w} \car_{v, w, >0}$, where for any $v, w
\in W$ such that $v \le w$,  $\car_{v, w, >0}$ is a semi-algebraic
cell and is a connected component of $\car_{v, w}$.

The following result generalizes these facts.

\proclaim{Corollary 3.16} $\overline{G_{>0}}=\bigsqcup\limits_{J
\subset I} \bigsqcup\limits_{\scriptstyle v, w, v', w' \in W^J
\atop \scriptstyle v \le w, v' \le w'} \bigsqcup\limits_{y, y' \in
W_J} Z^{v, w, v', w'; y, y'}_{J, >0}$. Moreover, for any $v, w,
v', w' \in W^J, y, y' \in W_J$ with $v \le w$, $v' \le w'$, $Z^{v,
w, v', w'; y, y'}_{J, >0}$ is a connected component of $Z^{v, w,
v', w'; y, y'}_J$ and is a semi-algebraic cell which is isomorphic
to $\bold R^d_{>0}$, where $d=l(w)+l(w')+2 l(w^J_0)+\mid J
\mid-l(v)-l(v')-l(y)-l(y')$.
\endproclaim

Proof. $\cp^J_{v, w, >0}$ (resp. $\cp^J_{v', w', >0}$) is a
connected component of $\cp^J_{v, w}$ (resp. $\cp^J_{v', w'}$)
(see [L3]). Thus $\{(P, Q, \gamma) \in Z^{v, w, v', w'; y, y'}_J
\mid P \in \cp^J_{v, w, >0}, \psi(Q) \in \cp^J_{v', w',
>0}\}$ is open and closed in $Z^{v, w, v', w'; y, y'}_J$. To prove
that $Z^{v, w, v', w'; y, y'}_{J, >0}$ is a connected component of
$Z^{v, w, v', w'; y, y'}_J$, it is enough to prove that $Z^{v, w,
v', w'; y, y'}_{J, >0}$ is a connected component of $\{(P, Q,
\gamma) \in Z^{v, w, v', w'; y, y'}_J \mid P \in \cp^J_{v, w,
>0}, \psi(Q) \in \cp^J_{v', w', >0}\}$.

Assume that $g \in G_{\hbox{\bf v}_+, \hbox{\bf w}, >0}, g' \in
G_{\hbox{\bf v}'_+, \hbox{\bf w}', >0}$ and $l \in L$. We have
that $(^g P_J)^{B^+}$ is the unique element $B \in \car_{v, w}$
that is contained in $^g P_J$(see 1.4). Therefore ${(^g
P_J)^{B^+}=^g B^+}$. Similarly, $(^g P_J)^{B^-}=^{g \dot{w}^J_0}
B^+$, $(^{\psi({g'} \i)} Q_J)^{B^+}=^{\psi({g'} \i) \dot{w}^J_0}
B^-$ and $(^{\psi({g'} \i)} Q_J)^{B^-}=^{\psi(g') \i} B^-$. Thus
$\po \Bigl((^g P_J)^{B^+}, ^{g l \psi(g')} \bigl((^{\psi({g'} \i)}
Q_J)^{B^+} \bigr) \Bigr)=\po(B^+, ^{l \dot{w}^J_0} B^-)$ and $\po
\Bigl((^g P_J)^{B^-}, ^{g l \psi(g')} \bigl((^{\psi({g'} \i)}
Q_J)^{B^-} \bigr) \Bigr)=\po(^{\dot{w}^J_0} B^+, ^l B^-)$.
Therefore we have that $(^g P_J, ^{\psi(g') \i} Q_J, g H_{P_J} l
U_{Q_J} \psi(g')) \in Z^{v, w, v', w'; y, y'}_J$ if and only if $l
\in B^+ \dot{y} \dot{w_0} B^+ \dot{w_0} \dot{w}^J_0 \cap
\dot{w}^J_0 B^+ \dot{y'} \dot{w_0} B^+ \dot{w_0}=B^+ \dot{y} B^-
\dot{w}^J_0 \cap \dot{w}^J_0 B^+ \dot{y'} B^-$.

Note that $L \cap B^+ \subset ^{\dot{w}^J_0} B^-$. Thus for any $x
\in W_J$, $(L \cap B^+) \dot{x} (L \cap B^+) \subset B^+ \dot{x}
\dot{w}^J_0 B^- \dot{w}^J_0$. Therefore, $$\eqalignno{L \cap B^+
\dot{y} B^- \dot{w}^J_0 &=\bigsqcup\limits_{x \in W_J} (L \cap
B^+) \dot{x} (L \cap B^+) \cap B^+ \dot{y} B^- \dot{w}^J_0 \cr
&=(L \cap B^+) \dot{y} \dot{w}^J_0 (L \cap B^+). \cr}$$

Similarly, $L \cap \dot{w}^J_0 B^+ \dot{y'} B^-=(L \cap B^-)
\dot{w}^J_0 \dot{y'} (L \cap B^-)$.

Then $\{(P, Q, \gamma) \in Z^{v, w, v', w'; y, y'}_J \mid P \in
\cp^J_{v, w, >0}, \psi(Q) \in \cp^J_{v', w',
>0}\}$ is isomorphic to $G_{v, w, >0} \times G_{v', w', >0} \times
\bigl( (L \cap B^+) \dot{y} \dot{w}^J_0 (L \cap B^+) \cap (L \cap
B^-) \dot{w}^J_0 \dot{y'} (L \cap B^-) \bigr)/Z(L)$. Note that
$\bigl( (L \cap B^+) \dot{y} \dot{w}^J_0 (L \cap B^+) \cap (L \cap
B^-) \dot{w}^J_0 \dot{y'} (L \cap B^-) \bigr) \cap L_{\ge
0}=U^-_{y w^J_0, >0} T_{>0} U^+_{w^J_0 y', >0}$. Therefore
$$\eqalignno{Z^{v, w, v', w'; y, y'}_{J, >0} & \cong G_{v, w,
>0} \times G_{v', w', >0} \times U^-_{y w^J_0, >0} T_{>0}
U^+_{w^J_0 y', >0}/\bigl(Z(L) \cap T_{>0} \bigr)\cr & \cong \bold
R^{l(w)+l(w')+2 l(w^J_0)+\mid J \mid-l(v)-l(v')-l(y)-l(y')}_{>0}.
\cr}$$

By 3.15, we have that $U^-_{y w^J_0, >0} T_{>0} U^+_{w^J_0 y',
>0}/\bigl(Z(L) \cap T_{>0} \bigr)$ is a connected component of
$\bigl( (L \cap B^+) \dot{y} \dot{w}^J_0 (L \cap B^+) \cap (L \cap
B^-) \dot{w}^J_0 \dot{y'} (L \cap B^-) \bigr)/Z(L)$. The corollary
is proved. \qed

\head Acknowledgements. \endhead I thank George Lusztig for
suggesting the problem and for many helpful discussions. I also
thank the referee for pointing out several mistakes in the
original manuscript and for some useful comments, especially
concerning 3.8, 3.10 and 3.15.

\enddocument